\nonstopmode
\documentclass[DIV=13]{scrartcl}
\usepackage[sb]{libertine} %
\usepackage[T1]{fontenc}
\usepackage{textcomp}
\usepackage[varqu,varl]{zi4}%
\usepackage[amsthm,upint]{libertinust1math} %
\usepackage[scr=boondoxo,bb=boondox]{mathalpha} %
\usepackage{mathtools}
\usepackage{todonotes}
\usepackage{booktabs}
\usepackage{enumitem}
\usepackage{placeins}
\usepackage{colortbl}
\usepackage{tikz}
\usepackage[citestyle=numeric-comp,minbibnames=5,maxbibnames=5,defernumbers=true,giveninits=true]{biblatex}
\usepackage{scrhack}
\usetikzlibrary{matrix}
\tikzset{sparsemat/.style={matrix of math nodes,execute at empty
cell={\node[gray]{0};},%
every left delimiter/.style={xshift=1ex},%
every right delimiter/.style={xshift=-1ex}, left delimiter={(},right
delimiter={)} } } %
\newenvironment{proofp}{\paragraph*{Proof.}}{\hfill $\blacksquare$ \\}
\usepackage{bm}
\newcommand{\symbf}[1]{\bm{#1}}
\newcommand{\symcal}[1]{\mathcal{#1}}
\newcommand{\symbfcal}[1]{\boldsymbol{\mathcal #1}}
\AtEveryBibitem{%
  \clearfield{urlyear}%
  \clearfield{urldate}%
}%
\AtEveryCitekey{%
  \clearfield{urlyear}%
  \clearfield{urldate}%
}%
\usepackage{xcolor}
\usepackage{float}
\RequirePackage{siunitx}
\RequirePackage[algoruled,norelsize, linesnumbered, lined, commentsnumbered]{algorithm2e}

\usepackage{mycommands}
\renewcommand{\mat}[1]{\mathbf{\uppercase{#1}}} %

\renewcommand{\vct}[1]{\mathbf{\lowercase{#1}}} %
\usepackage[format=plain,labelfont=bf]{caption}

\usepackage[pdftex,colorlinks=true,linkcolor={blue!50!black},citecolor={red!50!black},urlcolor=blue]{hyperref}

\renewcommand{\vec}[1]{\mathbf{#1}}

\definecolor{myblue}{RGB}{0 83 139}
\definecolor{myred}{RGB}{114 16 69}
\definecolor{mygreen}{RGB}{0 94 0}

\RequirePackage{tcolorbox}
\tcbuselibrary{skins,theorems,breakable}
\tcbsetforeverylayer{shield externalize}
\tcbset{
    thmstyle/.style={theorem style=standard,
      enhanced,breakable,
      arc=0mm,
      outer arc=0mm,
      boxrule=0mm,
      toprule=.7mm,
      bottomrule=.7mm,
      left=1mm,
      right=1mm,
      titlerule=0mm,
      toptitle=0mm,
      bottomtitle=1mm,
      top=0mm,
      colframe=black!9!white,
      colback=black!3!white,
      coltitle=black,
      title style={
        top color=black!5!white,
        middle color=black!5!white,
        bottom color=black!5!white},
      fonttitle=\bfseries\sffamily\normalsize,fontupper=\normalsize,},
    othstyle/.style={theorem style=standard,
      enhanced,breakable,
      arc=0mm,
      outer arc=0mm,
      boxrule=0mm,
      toprule=0mm,
      bottomrule=0mm,
      left=0mm,
      right=0mm,
      titlerule=0mm,
      toptitle=0mm,
      bottomtitle=1mm,
      top=0mm,
      colframe=white,
      colback=white,
      coltitle=black,
      title style={top color=white,bottom color=white,middle
        color=white},
      fonttitle=\bfseries\sffamily\normalsize,fontupper=\normalsize,}
  }
\tcolorboxenvironment{proof}{%
blanker,breakable,left=0mm,
before skip=10pt,after skip=10pt,}
  \newtcbtheorem[auto counter,number within=section]{definition}{Definition}{thmstyle}{def}
  \newtcbtheorem[auto counter,number within=section]{theorem}{Theorem}{thmstyle}{thm}
  \newtcbtheorem[auto counter,number within=section]{postulate}{Postulate}{thmstyle}{pos}
  \newtcbtheorem[auto counter,number within=section]{axiom}{Axiom}{thmstyle}{axi}
  \newtcbtheorem[auto counter,number within=section]{corollary}{Corollary}{thmstyle}{cor}
  \newtcbtheorem[auto counter,number within=section]{lemma}{Lemma}{thmstyle}{lem}
  \newtcbtheorem[auto counter,number within=section]{proposition}{Proposition}{thmstyle}{pro}
  \newtcbtheorem[auto counter,number within=section]{conclusion}{Conclusion}{thmstyle}{conc}
  \newtcbtheorem[auto counter,number within=section]{remark}{Remark}{thmstyle}{rem}
  \newtcbtheorem[auto counter,number within=section]{example}{Example}{othstyle}{ex}
  \newtcbtheorem[auto counter,number within=section]{problem}{Problem}{thmstyle}{problem}
  \newtcbtheorem[auto counter,number within=section]{assumption}{Assumption}{thmstyle}{ass}
  \newtcbtheorem[auto counter,number within=section]{notation}{Notation}{thmstyle}{not}
  \newtcbtheorem[auto counter,number within=section]{algo}{Algorithm}{thmstyle}{alg}

  \date{}
\begin{document}
\title{Well-posedness and exponential stability of dispersive nonlinear Maxwell equations with PML: An evolutionary approach}
\author{Nils Margenberg
  \thanks{Helmut Schmidt University, Faculty of Mechanical and Civil
    Engineering, Holstenhofweg 85, 22043 Hamburg, Germany} \thanks{Corresponding
    author: \href{mailto:margenbn@hsu-hh.de}{\texttt{margenbn@hsu-hh.de}}}%
  \and Markus
  Bause\footnotemark[1] }
\maketitle
\vspace{-7ex}

\begin{abstract}
  This paper presents a mathematical foundation for physical models in nonlinear
  optics through the lens of evolutionary equations. It focuses on two key
  concepts: \emph{well-posedness} and \emph{exponential stability} of Maxwell
  equations, with models that include materials with complex dielectric
  properties, dispersion, and discontinuities. We use a Hilbert space framework
  to address these complex physical models in nonlinear optics. While our focus
  is on the first-order formulation in space and time, higher solution
  regularity recovers and equates to the second-order formulation. We
  incorporate perfectly matched layers (PMLs), which model absorbing boundary
  conditions, to facilitate the development of numerical methods. We demonstrate
  that the combined system remains well-posed and exponentially stable. Our
  approach applies to a broad class of partial differential equations (PDEs) and
  accommodates materials with nonlocal behavior in space and time. The
  contribution of this work is a unified framework for analyzing wave
  interactions in advanced optical materials.

\end{abstract}

\section{Introduction}
The accurate modeling of wave interactions with dispersive materials is
essential in the development of photonic devices, optical waveguides,
electromagnetic compatibility, and metamaterials. Nonlinear optics, in
particular, deals with phenomena where the response of a material to
electromagnetic fields is nonlinear, leading to effects such as harmonic
generation, self-focusing, and soliton propagation. A significant challenge in
this field is the modeling of wave propagation in media with complex
permittivity and permeability, as well as the characterization of the wave
phenomena that occur. In this paper, we present a mathematical foundation for
modeling such complex physical systems by focusing on the well-posedness and
exponential stability of nonlinear, dispersive electromagnetic wave equations
within the framework of evolutionary equations. We utilize a Hilbert space
framework and adopt a first-order formulation in space and time. We show
well-posedness and stability for a general class of material laws, which
incorporate memory effects and nonlocal responses in material models.

In particular, we employ the framework of
evolutionary equations developed by R.
Picard~\autocite{picardStructuralObservationLinear2009}. 
We consider physical systems, which can be modeled by the abstract PDE
\begin{equation}
\label{eq:evolutionary-problem}
(\partial_t \symbfcal{M}_0 + \symbfcal{M}_1 + \symbfcal{A})\vct{u} = \vct{f}\,,
\end{equation}
where \(\symbfcal{A}\) is a skew self-adjoint operator in a Hilbert space
\(\symbfcal{H}\), and \(\symbfcal{M}_0\) and \(\symbfcal{M}_1\) are linear and
bounded operators in \(\symbfcal{H}\). This abstract framework allows us to
address the well-posedness and exponential stability of nonlinear, dispersive
Maxwell equations. The theory is based exponentially in time weighted function
spaces on the entire axis, thereby avoiding the necessity of an initial
condition. The abstract PDE~\eqref{eq:evolutionary-problem} accommodates a general class of
material laws, including those with complex dielectric properties, anisotropic
and discontinuous material coefficients, as well as fields with discontinuities.
Moreover, it enables the analysis of materials exhibiting nonlocal behavior in
both space and time, which is essential for accurately modeling advanced optical
materials like metamaterials and engineered media.

In order to truncate unbounded domains to bounded domains, we incorporate PMLs
to model absorbing boundary conditions. PMLs enable the truncation to bounded
domains without reflecting waves back into the interior, which is crucial for
numerical simulations. We address the additional complexity arising from the
coupling of physical and PML domains and demonstrate that the combined system is
well-posed and exponentially stable. Traditional methods, such as semigroup
theory or Galerkin approximations, exhibit difficulties to prove well-posedness
for these coupled models. An advantage of the theory of R. Picard is its
reliance on mild assumptions for the involved operators, making it suitable for
effectively handling coupled physical and PML models. The theory developed by R.
Picard has the advantage of requiring only mild assumptions for the involved
operators, making it an attractive option for handling coupled physical and PML
models.

The central result, Theorem~\ref{thm:picards}, first introduced by
Picard~\autocite{picardStructuralObservationLinear2009}, applies to a broad
class of PDEs modeling various physical phenomena. This generality is crucial
for addressing systems of changing type and their
discretization~\autocite{franzContinuousTimeIntegration2021,
  franzNumericalMethodsChanging2019}. We extend the existing well-posedness
results to a broad class of practical dispersion models and their integration
with PMLs. While the combination of PML and nonlinear models has been used in
practice~\autocite{margenbergAccurateSimulationTHz2023,
  margenbergOptimalDirichletBoundary2024,
  abrahamConvolutionFreeMixedFiniteElement2019}, these formulations have not
been thoroughly investigated in the context of theoretical analysis. The
formulations presented here are suitable for the development of numerical
discretizations and serve as an initial step towards bridging the gap between
mathematical theory and numerical simulations in nonlinear optics.

This paper proceeds with a review of relevant literature in this section, after
which the essential notation is introduced in
Section~\ref{sec:preliminaries-notation}. In Section~\ref{sec:picard}, we
introduce the framework of evolutionary equations and present the key theorem
(Theorem~\ref{thm:picards}) underpinning our well-posedness results. We further
review the well-posedness of models in nonlinear optics, building upon recent
results by Ionescu-Tira and
Dohnal~\autocite{dohnalWellPosednessExponentialStability2023}. In
Section~\ref{sec:dispersion-picard} we introduce a class of dispersive Maxwell
equations, show their well-posedness and discuss their exponential stability. We
introduce a formulation of PML in the framework of evolutionary equations in
Section~\ref{sec:pml-picard} and demonstrate their well-posedness and
exponential stability for a range of dispersion models in
Section~\ref{sec:dispersion-pml-picard}.

\subsection*{Related Works}\label{sec:related}

\paragraph*{Well-Posedness and Error Analysis}
Well-posedness and error analysis for Maxwell's equations, particularly in nonlinear optics, is an active research area. Nonlinear effects arise from the interaction of light with material electrons, modeled as nonlinear polarization. Similar nonlinear models appear in acoustic and elastic wave equations, such as the Westervelt or Kuznetsov equations~\autocite{antoniettiHighorderDiscontinuousGalerkin2020}.
Quasilinear hyperbolic evolution equations with Kerr nonlinearity have been analyzed within the Kato framework~\autocite{katoQuasilinearEquationsEvolution1975}, with refinements under relaxed data assumptions~\autocite{mllerWellposednessGeneralClass2014}. Numerical error analysis using implicit Euler methods has been conducted for these equations~\autocite{hochbruckErrorAnalysisImplicit2017}.
Local well-posedness for Maxwell's equations in Kerr media with perfectly conducting boundary conditions has been established using energy techniques~\autocite{spitzLocalWellposednessNonlinear2019, spitzLocalWellposednessNonlinear2017}. Further developments include error analysis providing stability and error bounds in Hilbert spaces~\autocite{hochbruckErrorAnalysisSpace2022, kovcsStabilityConvergenceTime2018}. Optimal-order error bounds have been obtained for discretizations using isoparametric finite elements and various time discretizations~\autocite{drichStrongNormError2022, drichErrorAnalysisExponential2021}.
Within the abstract evolutionary equation framework, the well-posedness of general Maxwell's equations has been explored~\autocite{picardStructuralObservationLinear2009, seifertEvolutionaryEquationsPicard2022, dohnalWellPosednessExponentialStability2023}. Well-posedness and error estimates for the nonlinear Helmholtz equation in Kerr media using finite element methods have also been examined~\autocite{verfrthHigherorderFiniteElement2023, maierMultiscaleScatteringNonlinear2022}. The Kato framework has proven instrumental in establishing stability and error bounds, including first-order convergence of implicit Euler schemes for Maxwell's equations with Kerr-type nonlinearity~\autocite{drichStrongNormError2022}.
In the time-harmonic setting, well-posedness of the nonlinear electromagnetic wave equation in Kerr media has been shown, with error estimates provided for finite element methods~\autocite{verfrthHigherorderFiniteElement2023}. A multiscale approach to the nonlinear Helmholtz problem has demonstrated well-posedness and convergence results~\autocite{maierMultiscaleScatteringNonlinear2022}.

\paragraph*{Perfectly Matched Layers}
The concept of perfectly matched layers (PMLs), introduced by Bérenger~\autocite{berengerPerfectlyMatchedLayer1994}, significantly advanced numerical simulations of wave propagation. PMLs offer straightforward implementation compared to higher-order absorbing boundary conditions and perform well across various applications, including acoustics and elastodynamics. However, the mathematical analysis of PMLs remains challenging, with certain physical models leading to unstable formulations. The ongoing debate between PMLs and absorbing boundary conditions (ABCs) continues without a definitive conclusion.
PMLs can be implemented through field splitting and interface treatments or considered in the frequency domain as coordinate stretching~\autocite{chewComplexCoordinateStretching1997}, leading to time-dependent PMLs computable via auxiliary differential equations. For a comprehensive review of PMLs, refer to~\autocite{pledReviewRecentDevelopments2021}.

\paragraph*{Evolutionary Equations}
The seminal work by Picard~\autocite{picardStructuralObservationLinear2009} established the well-posedness of evolutionary equations such as~\eqref{eq:evolutionary-problem}, providing a foundation for a wide range of applications~\autocites{waurickStabilizationHomogenization2016, picardComprehensiveClassLinear2016, picardClassBoundaryControl2012, franzHomogenisationParabolicHyperbolic2020, picardModelsElasticSolids2015, picardEvolutionaryEquationsMaterial2015, picardElectromagneticWavesComplex2013} in both linear and nonlinear systems.
A key concept in this framework is exponential stability, where exponentially decaying right-hand sides lead to exponentially decaying solutions~\autocites{trostorffExponentialStabilityLinear2013a, trostorffExponentialStabilityInitial2018b}. This is particularly relevant for nonautonomous problems~\autocites{dohnalWellPosednessExponentialStability2023, picardWellposednessClassNonautonomous2018, sSolutionTheoryGeneral2017, trostorffWellposednessGeneralClass2020}.
Picard's framework has significantly influenced the numerical discretization of evolutionary equations~\autocites{franzPostprocessingImprovedError2023, franzNumericalMethodsChanging2019, franzContinuousTimeIntegration2021}, especially for systems of changing type~\autocite{franzNumericalMethodsChanging2019}. Formulating problems in spaces of higher spatial regularity accommodates a broader class of nonlinearities~\autocite{dohnalWellPosednessExponentialStability2023}, as demonstrated by research on spatial regularity for evolutionary equations~\autocites{picardMaximalRegularityClass2017, trostorffMaximalRegularityNonautonomous2021}.

Friedrichs proposed a similar approach for elliptic and hyperbolic problems,
formulating them as abstract operator equations
\(\symbfcal{K}\vct{u} = \vct{f}\) with an accretive symmetric operator
\(\symbfcal{K}\)~\autocites{friedrichsSymmetricPositiveLinear1958}{friedrichsSymmetricHyperbolicLinear1954}.
These systems have been extensively discussed, particularly for the parabolic
case~\autocites{antoniHeatEquationFriedrichs2013}{antoniSecondorderEquationsFriedrichs2014}{burazinNonStationaryAbstractFriedrichs2016}{antoniFriedrichsSystemsHilbert2017}{antoniComplexFriedrichsSystems2017},
but they often focus on local operators in space and neglect non-stationary
mixed-type examples. A significant limitation of Friedrichs systems is their
inability to distinguish between time and space coordinates, thus neglecting
causality in time evolution. In contrast, Picard's approach inherently
incorporates causality due to a uniform positive definiteness constraint on the
operators~\autocites{picardStructuralObservationLinear2009}{seifertEvolutionaryEquationsPicard2022}.
Inspired by Friedrichs systems, comprehensive numerical analyzes of various
PDEs, including changing-type equations, have been
performed~\autocites{ernDiscontinuousGalerkinMethods2006}{ernDiscontinuousGalerkinMethods2006c}{ernDiscontinuousGalerkinMethods2008}{jensenDiscontinuousGalerkinMethods2004}.

\section{Preliminaries and notation}
\label{sec:preliminaries-notation}
Throughout this paper, we utilize several notations, which we introduce here.
Let $\Omega \subseteq \mathbb{R}^d$ be a bounded domain with
Lipschitz-continuous boundary, and let $\symbfcal{H}$ denote a real Hilbert
space. We use standard notations: $\dom(\symbfcal{A})$ denotes the domain of an
operator $\symbfcal{A}$ and $\symbfcal{B}(\symbfcal{H})$ denotes the space of
bounded linear operators on $\symbfcal{H}$.

We denote by $L^2(\Omega)$ the space of square-integrable functions over $\Omega$, and by $L^2(\Omega)^d$ the space of $d$-tuples of such functions. The inner product in $\symbfcal{H}$ is denoted by $\langle \cdot, \cdot \rangle_{\symbfcal{H}}$, and the associated norm by $\| \cdot \|_{\symbfcal{H}}$.
We employ the Fourier-Laplace transform $\symcal{L}_{\nu}$ with exponential weight $\nu$, defined by
\begin{equation}
  \label{eq:flapt}
(\symcal{L}_{\nu} \vct f)(\omega) = \frac{1}{\sqrt{2\pi}} \int_{\mathbb{R}} \vct f(t) \exp(-(\iu \omega + \nu) t) \, \drv t, \quad \vct f \in L^2_{\nu}(\mathbb{R}; \symbfcal{H}),
\end{equation}
where $\iu$ denotes the imaginary unit satisfying $\iu^2 = -1$. We will
frequently refer to Lipschitz continuous functions. A function $\vct{f}$ mapping
between normed vector spaces is said to be \emph{Lipschitz continuous} if there
exists a constant $L \geq 0$ such that
\[
\| \vct{f}(\vct u) - \vct{f}(\vct v) \| \leq L \| \vct u - \vct v \|, \quad \text{for all } \vct u, \vct v \in \dom(\vct{f}).
\]
The smallest such constant $L$ is called the \emph{Lipschitz constant} of $\vct{f}$, denoted by $\| \vct{f} \|_{\mathrm{Lip}}$, and is defined by
\[
\| \vct{f} \|_{\mathrm{Lip}} = \inf \left\{ L \geq 0 \mid \| \vct{f}(\vct u) - \vct{f}(\vct v) \| \leq L \| \vct u - \vct v \|, \ \forall \vct u, \vct v \right\}.
\]

\section{A Hilbert space framework for evolutionary
  equations\label{sec:picard}}
The framework of evolutionary equations introduced
in~\autocite{picardStructuralObservationLinear2009} covers a wide range of
evolution problems of mathematical physics. We consider problems arising from
electromagnetic wave propagation. For \(t \in \R\), the abstract evolutionary
equation is of the form presented in~\eqref{eq:evolutionary-problem}. Recall
that \(\symbfcal A\) is skew self-adjoint, \(\symbfcal M_0\) and
\(\symbfcal M_1\) are linear and bounded. The skew self-adjointness of
\(\symbfcal A\) typically arises from the block operator form
\begin{equation}
\label{eq:A-block-matrix}
\symbfcal A=\begin{pmatrix}
       0 & -\mat C_0^{*} \\ \mat C_0 & 0
    \end{pmatrix}\,,
\end{equation}
with a closed linear operator \(\mat C\colon \mat X_0\to \mat Y\) between real Hilbert
spaces \(\mat X_0\) and \(\mat Y\), where the Hilbert space \(\mat X_1=\dom(\mat C)\) is equipped with
the graph inner product of \(\mat C\)
(cf.~\autocite{picardStructuralObservationLinear2009,seifertEvolutionaryEquationsPicard2022}).
The skew-selfadjointness of \(\symbfcal A\) is a consequence of the
definition of adjoints and the fact that \(\mat C = \mat C^{**}\). Therefore, a concrete
Hilbert space \(\symbfcal{H}\) we introduced above is defined as
\(\symbfcal{H}\coloneq \mat X_0\oplus \mat Y\). Our aim in utilizing the framework of
evolutionary equations is two-fold. First, we want to prove the well-posedness
of dispersive Maxwell equations relevant in the field of optics. In follow up work,
we want to develop numerical methods for the arising first order formulations of
the aforementioned problems.

To this end, we use the solution theory
from~\autocite{picardStructuralObservationLinear2009} (see
also~\autocites[Theorem
6.2.1]{seifertEvolutionaryEquationsPicard2022}[][]{picardPartialDifferentialEquations2011}),
which allows to achieve these two goals in one holistic framework. An essential
ingredient for the well-posedness theory are exponentially weighted $L^2$
spaces.
\begin{definition}{Exponentially weighted \(L^2\) space}{exp-weight-l2}
  Let \(\symbfcal{H}\) be a Hilbert space over \(\R\). For \(\nu\in \R\), \(\nu > 0\) we
  consider a weighting function \(t\mapsto \exp(-\nu t)\) to define the Hilbert
  space
  \begin{equation}
    \label{eq:exp-wgt-hilbert}
    L^2_{\nu}(\R;\: \symbfcal{H})\coloneq \brc{\vct f\colon \R \to \symbfcal{H}
      \suchthat \vct f\text{ is measurable},\;\int_{\R} \norm{\vct f(t)}_{\symbfcal{H}}^{2} \exp(-2\nu t)\drv t < \infty}.
  \end{equation}
  The space \(L^2_{\nu}(\R;\:\symbfcal{H})\) is a Hilbert space with respect to the inner
  product
  \begin{equation}
    \label{eq:wgt-hilbert-inner}
    \dd{\vct f}{\vct g}_{\nu}\coloneq\int_{\R}\dd{\vct f(t)}{\vct g(t)}_{\symbfcal{H}}\exp(-2\nu t)\drv t\,,\quad\forall \vct f,\: \vct g\in L^2_{\nu}(\R;\:\symbfcal{H}).
  \end{equation}
  We denote the norm induced by the inner product~\eqref{eq:wgt-hilbert-inner}
  by \(\norm{\cdot}_{\nu}\).
\end{definition}
  \begin{definition}{Time derivative in exponentially weighted \(L^2\) space}{exp-weight-tdl2}
  We define the derivative with respect to time as the closure of the operator
  \begin{equation}
    \label{eq:partial-t}
    \begin{aligned}
      \partial_t\colon C^{\infty}_{c}(\R;\: \symbfcal{H}) \subseteq L^2_{\nu}(\R;\: \symbfcal{H})&\to L^2_{\nu}(\R;\: \symbfcal{H})\\
      \vct u&\mapsto \vct u'\,,
    \end{aligned}
  \end{equation}
  where \(C^{\infty}_{c}(\R;\: \symbfcal{H})\) is the space of infinitely differentiable
  \(\symbfcal{H}\)-valued functions on \(\R\) with compact support (denoted by the
  subscript $c$). The domain of \(\partial_t^s\) for \(s\in \N\) is denoted by
  \(H^s_{\nu}(\R;\:\symbfcal{H})\).
\end{definition}
\begin{definition}{Curl operator}{mw-operator}
  Let \(\Omega\subseteq \R^3\), be a bounded domain with
  Lipschitz-continuous boundary. We introduce
  $\curl_0$ as the closure of the operator $\curl_c$ in $L^2(\Omega)^3$, which is defined by
\begin{align}
\curl_c\colon C_c^\infty(\Omega)^{3}\subseteq L^2(\Omega)^{3}  &\to L^2(\Omega)^{3}\\
\symbf\varphi=\begin{pmatrix}
    \varphi_x \\ \varphi_y \\ \varphi_z
  \end{pmatrix}
  &\mapsto
  \nabla\times\symbf\varphi =
  \begin{pmatrix}
    \partial_y \varphi_{z} - \partial_z\varphi_y \\
    \partial_z \varphi_{x} - \partial_x\varphi_z \\
    \partial_x \varphi_{y} - \partial_y\varphi_x
  \end{pmatrix}\,.
\end{align}
As in~\autocite[Section 6.1]{seifertEvolutionaryEquationsPicard2022} and~\autocite[Section 2.3]{giraultFiniteElementMethods1986}, we use
\begin{equation*}
  \dom(\curl_0) = \mat H_0(\curl,\,\Omega) \coloneq
  \overline{C_c^\infty(\Omega)^3}^{\lVert \cdot \rVert_{\mat H(\curl)}}\,,
\end{equation*}
with the norm $\lVert \boldsymbol \varphi \rVert_{\mat H(\curl,\,\Omega)} = \bigl( \lVert \boldsymbol \varphi
\rVert_{L^2(\Omega)}^2 + \lVert \nabla\times \boldsymbol \varphi \rVert_{L^2(\Omega)}^2
\bigr)^{1/2}$.
We further set $\curl \coloneq \curl_0^*$ and introduce the notation
\begin{equation*}
  \dom(\curl) = \mat H(\curl,\,\Omega) \coloneq \{ \vct u \in L^2(\Omega)^3 : \curl \vct u \in L^2(\Omega)^3 \}\,.
\end{equation*}
Note that $\mat H_0(\curl,\,\Omega)\subseteq \mat H(\curl,\,\Omega)$ is the space of
$\mat H(\curl,\,\Omega)$ functions with vanishing tangential component
(cf.~\autocite[Remark 6.1.3]{seifertEvolutionaryEquationsPicard2022}).
\end{definition}
Given the \(\curl\) and \(\curl_0\) operator, we consider the following operator
throughout this paper
\begin{equation}\label{eq:curl-A}
  \symbfcal A \coloneq
  \begin{pmatrix}
    {\color{gray} \symbf 0} & -\curl \\ \curl_0 & {\color{gray} \symbf 0}
  \end{pmatrix}.
\end{equation}
Then, since \(\curl^* = (\curl_0^*)^* = \curl_0\), \(\symbfcal A \colon
\mat H_0(\curl,\,\Omega)\times \mat H(\curl,\,\Omega) \subseteq L^2(\Omega)^3 \to
L^2(\Omega)^3\) is skew self-adjoint, i.\,e.\
\begin{equation*}
  \symbfcal A^* =
  \begin{pmatrix}
    {\color{gray} \symbf 0} & -\curl \\ \curl_0 & {\color{gray} \symbf 0}
  \end{pmatrix}^* =
  \begin{pmatrix}
    {\color{gray} \symbf 0} & \curl_0^* \\ -\curl^* & {\color{gray} \symbf 0}
  \end{pmatrix} =
  \begin{pmatrix}
    {\color{gray} \symbf 0} & \curl \\ -\curl_0 & {\color{gray} \symbf 0}
  \end{pmatrix} = -\symbfcal A\,.
\end{equation*}
The above choice of
\begin{equation*}
  \dom(\symbfcal A) = \mat H_0(\curl,\,\Omega)\times \mat H(\curl,\,\Omega)
\end{equation*}
and the skew self-adjointness of \(\symbfcal A\) already encode eventually
appearing interface and boundary conditions of a perfect conductor
(cf.~\autocite{dohnalWellPosednessExponentialStability2023}). Throughout this
work, we consider a Cauchy problem in \(L^2(\Omega)^3\times L^2(\Omega)^3\) for the Maxwell system
with $\vct u = (\vct e,\,\vct h)$ formulated as an evolutionary problem
in positive time:
\begin{equation}
    \begin{aligned}
      \partial_t
      \begin{pmatrix}
        \vct D(\vct E) \\ \vct  B(\vct H)
      \end{pmatrix} +
      \begin{pmatrix}
        {\color{gray} \symbf 0} & -\curl \\ \curl_0 & {\color{gray} \symbf 0}
      \end{pmatrix}
      \begin{pmatrix}
        \vct E \\ \vct H
      \end{pmatrix} &=
      \begin{pmatrix}
        -\vct J \\ {\color{gray} \symbf 0}
      \end{pmatrix}\,, & t > 0\,,\\
      (\vct E(t),\,\vct H(t)) &= (\vct E_0(t),\,\vct H_0(t)), & t \le 0\,.
    \end{aligned}
  \label{eq:Maxwell_pos_init}
\end{equation}
Here, let
\begin{equation*}
(\vct E_0,\,\vct H_0) \colon (-\infty,\,0] \to L^2(\Omega)^3\times L^2(\Omega)^3,\quad t\mapsto (\vct E_0(t),\,\vct H_0(t))\,,
\end{equation*}
denote a given history. This is indispensable as the material function
\((\vct{D},\,\vct{B})\) often depends on past values of its arguments. To
satisfy jump conditions for \(\vct{E}\) and \(\vct{H}\),
\((\vct{E}(t),\,\vct{H}(t))\) must remain within the domain of \(\symbfcal{A}\)
for all \(t\) (cf.~\autocite{dohnalWellPosednessExponentialStability2023}).
As discussed in~\autocite{dohnalWellPosednessExponentialStability2023},
the interface setting does not impact the solution theory, since transmission
conditions are inherently incorporated within the domain \(\dom(\symbfcal{A})\).

Note that the divergence equations of \(\vct{D}\) and \(\vct{B}\) are
redundant, as they can be derived from~\eqref{eq:Maxwell_pos_init}
given appropriate initial conditions.
Spitz~\autocite{spitzLocalWellposednessNonlinear2019} elaborates that by
applying the divergence operator to the first term
in~\eqref{eq:Maxwell_pos_init} and integrating, we find \(\Dv \vct{D}(t) =
\varrho(t)\) holds for \(t \ge 0\) if, and only if, \(\varrho\) and \(\vct{J}\)
satisfy
\begin{equation*}
\varrho(t) = \varrho(0) - \int_0^t \Dv \vct{J}(s) \, \drv s\,.
\end{equation*}
Similarly, from the second term in~\eqref{eq:Maxwell_pos_init}, \(\Dv \vct{B}\)
remains constant for all \(t > 0\). Therefore, if \(\Dv \vct{B}(0) = 0\), then
\(\Dv \vct{B}(t) = 0\) holds for all \(t > 0\).

\subsection{Well-posedness of linear evolutionary problems \label{ssec:wp-evolutionary-problem}}
\label{sec:orgfe4224b}
In~\eqref{eq:evolutionary-problem} we introduced \(\symbfcal{A}\) as a
skew-selfadjoint operator in a Hilbert space \(\symbfcal{H}\), and
\(\symbfcal{M}_0\) and \(\symbfcal{M}_1\) as linear and bounded operators in
\(\symbfcal{H}\). We refine this by defining a \emph{linear material law}
\(\symbfcal{M}\), where
\(\symbfcal{M}_0 + z^{-1}\symbfcal{M}_1 \coloneqq \symbfcal{M}\). In the context
of the Maxwell equation~\eqref{eq:Maxwell_pos_init} the material law
$\symbfcal{M}$ encapsulates the functions $\vct D(\vct E)$ and $\vct B(\vct H)$.
\begin{definition}{Linear material law~\autocite[Section 5.3]{seifertEvolutionaryEquationsPicard2022}}{matlaw}
  A \emph{linear material law} \(\symbfcal{M}\) is an analytical mapping, with
  \(\symbfcal{M} \colon \dom(\symbfcal{M}) \subseteq \mathbb{C} \to
  \symbfcal{B}(\symbfcal{H})\). The material law \(\symbfcal{M}\)
  is uniformly bounded within a right half-plane:
\begin{equation*}
  \exists \nu_0 \in \mathbb{R} \colon \sup_{\re z > \nu_0} \lVert
  \symbfcal{M}(z) \rVert_{\symbfcal{H}}\eqcolon \lVert \symbfcal{M}(z)
  \rVert_{\infty,\C_{\re z > \nu}}< \infty\,.
\end{equation*}

The operator \(\symbfcal{M}(\partial_t)\) is defined by:
\begin{equation*}
\symbfcal{M}(\partial_t) = \symcal{L}_{\nu}^* \symbfcal{M}(\iu \cdot + \nu) \symcal{L}_\nu,
\end{equation*}
where \(\symbfcal{M}(\iu \cdot + \nu)\) is
\((\symbfcal{M}(\iu \cdot + \nu) \boldsymbol \varphi)(t) = \symbfcal{M}(\iu t + \nu)
\boldsymbol \varphi(t)\), and
\(\symcal{L}_{\nu} \colon L^2_{\nu}(\mathbb{R},\,\symcal{H}) \to
L(\mathbb{R},\,\symcal{H})\) is the unitary extension of the Fourier-Laplace
transform~\eqref{eq:flapt}.
\end{definition}
We identify \(\symbfcal{M}_0\), \(\symbfcal{M}_1\), and \(\symbfcal{A}\) with
their canonical extensions to \(\symbfcal{H}\)-valued functions acting as
abstract multiplication operators. We now introduce the concept of
\emph{causality}, which is essential for the solution theory of evolutionary
equations.
\begin{definition}{Causality}{causa}
Given \( a \in \R \), the multiplication operator \( \theta_a^+ \colon L^2_{\nu}(\R, \symbfcal{H}) \to L^2_{\nu}(\R, \symbfcal{H}) \) is defined as
\[
\theta_a^+ \vct u(t) = \chi_{(a,\,\infty)}(t) \vct u(t) =
\begin{cases}
\vct u(t), & \text{if } t > a, \\
0, & \text{if } t \leq a,
\end{cases}
\]
where \( \chi_{(a,\,\infty)} \) is the characteristic function of the interval \( (a,\,\infty) \).
A mapping \(\vct f \in L^2_{\nu}(\R, \symbfcal{H}) \) is \emph{(forward) causal} if, for all \( a \in \R \),
\[
(1 - \theta_a^+)(\vct u - \vct v) = 0 \implies (1 - \theta_a^+)(\vct f(\vct u) - \vct f(\vct v)) = 0.
\]
That is, if \( \vct u \) and \(\vct v \) coincide on \( (-\infty,a] \), then \( \vct f(\vct u) \) and \( \vct f(\vct v) \) also coincide on \( (-\infty,a] \).
\end{definition}
Given a linear material law \(\symbfcal{M}\), the operator
\(\symbfcal{M}(\partial_t)\) is causal on \(L^2_\nu(\R, \symbfcal{H})\) for
\(\nu > \nu_0\) due to the Paley-Wiener theorem~\autocite[Theorem
8.1.2]{seifertEvolutionaryEquationsPicard2022}. The operator \(\partial_t\) is
boundedly invertible for \(\nu \ne 0\) and causally invertible for \(\nu > 0\).
For \(\nu > 0\), \(\partial_t^{-1} \colon L^2_\nu(\R, \symbfcal{H}) \to
L^2_\nu(\R, \symbfcal{H})\) is given by
\begin{equation}
(\partial_t^{-1} \vct u)(t) = \int_{-\infty}^t \vct u(s) \, \drv s\,,
\label{eq:Antiderivative}
\end{equation}
with \(\lVert \partial_t^{-1} \rVert \le 1/\nu\)~\autocite[Section 3.2]{seifertEvolutionaryEquationsPicard2022}. The notation \(\partial_t^{-1}\) is reserved for the causal map defined in~\eqref{eq:Antiderivative}.
\begin{theorem}{Well-posedness of evolutionary equations~\autocite[Theorem 6.2.1]{seifertEvolutionaryEquationsPicard2022}}{picards}
  Let $\symbfcal A: \dom(\symbfcal A) \subseteq \symbfcal{H} \to \symbfcal{H}$ be a skew self-adjoint operator, and let $\symbfcal M$ represent a linear material law. Suppose $z\symbfcal M(z)=z\symbfcal M_0(z)+\symbfcal M_1(z)$ is strictly accretive on a half-plane $\C_{\re > \nu_0}$ for some $\nu_0 \in \R$, such that
\begin{equation}
  \exists \gamma > 0 \ \forall z \in \C_{\re > \nu_0}: \ \re z\symbfcal M(z) \ge \gamma,
  \label{eq:MaterialLawAccretivity}
\end{equation}
or, equivalently, $\re \langle z\symbfcal M(z)\vct u, \vct u \rangle \ge \gamma\lVert \vct u \rVert_{\symbfcal{H}}^2$ for all $\vct u \in \symbfcal{H}$.
Then, for any $\nu > \nu_0$, the operator ${\partial_t\symbfcal M(\partial_t) + \symbfcal A}$ is closable, and the mapping
\begin{equation}
  \symbfcal{S}_\nu \coloneq {(\overline{\partial_t\symbfcal M(\partial_t) +
      \symbfcal A})}^{-1} \colon L^2_\nu(\R,\,\symbfcal{H}) \to
  L^2_\nu(\R,\,\symbfcal{H})
    \label{eq:inv-op}
\end{equation}
is well-defined, bounded, and satisfies $\lVert \symbfcal{S}_\nu \rVert_{L^2_\nu\to
  L^2_\nu} \le\frac{1}{\gamma} $. Here the closure is taken in $L^2_\nu(\R;\: \symbfcal{H})$.
Furthermore, $\symbfcal{S}_\nu$ is causal and for every $\vct f \in L^2_\nu(\R,\,\symbfcal{H})\cap L^2_{\nu'}(\R,\,\symbfcal{H})$, $\symbfcal{S}_\nu \vct f = \symbfcal{S}_{\nu'}\vct f \in L^2_\nu(\R,\,\symbfcal{H})\cap L^2_{\nu'}(\R,\,\symbfcal{H})$ whenever $\nu,\,\nu' > \nu_0$.
\end{theorem}
We remark the following on Theorem~\ref{thm:picards}:
\begin{itemize}\itemsep1pt \parskip0pt \parsep0pt
\item The mapping $\symbfcal{S}_\nu$ in~\eqref{eq:inv-op} can be expressed in
  terms of the spectral representation of the time-derivative:
  \[
    \vct u = {(\overline{\partial_t\symbfcal M(\partial_t) + \symbfcal A})}^{-1}\vct f = \symcal{L}_\nu^{-1}\left((i\,\cdot + \nu) \symbfcal M(i\,\cdot + \nu) + \symbfcal A\right)^{-1}\symcal{L}_\nu \vct f,
  \]
  where \(\vct f \in L^2_\nu(\R,\symbfcal{H})\) and \(\nu > \nu_0\). This establishes sufficient conditions for the operator
  \[
    (\cdot\,\symbfcal M(\cdot) + \symbfcal A)^{-1} \colon \C_{\re > \nu_0} \cap \dom(\symbfcal M) \to \symbfcal{B}(\symbfcal{H}),
  \]
  to have a bounded and analytic extension in \(\C_{\re > \nu_0}\). Thus,~\eqref{eq:evolutionary-problem} is called \emph{well-posed} in \(L^2_\nu(\R,\symbfcal{H})\).
\item The stability estimate $\displaystyle \norm{\vct U}_\nu \leq \frac{1}{\gamma} \norm{\vct F}_\nu$ holds.
\item Let \(\symbfcal M_0, \symbfcal M_1 \colon \symbfcal{H} \to \symbfcal{H}\) be bounded linear operators, with \(\symbfcal M_0\) self-adjoint and positive definite, and \(\symbfcal A \colon \dom(\symbfcal A) \subseteq \symbfcal{H} \to \symbfcal{H}\) skew self-adjoint. The material law \(\symbfcal{M} \coloneq \symbfcal{M}_0 + z^{-1} \symbfcal{M}_1\) satisfies:
  \[
    \re \dd{\vct u}{z\symbfcal{M}(z) \vct u} = \dd{\vct u}{\nu \symbfcal M_0 \vct u} + \dd{\vct u}{\re \symbfcal M_1 \vct u}.
  \]
  Theorem~\ref{thm:picards} can be stated equivalently by replacing~\eqref{eq:MaterialLawAccretivity} with:
  \begin{equation}
    \label{eq:evo-cond}
    \exists \gamma > 0 \;\forall \nu \geq \nu_0, \;\vct u \in \symbfcal{H} \colon \langle (\nu \symbfcal M_0 + \re \symbfcal M_1) \vct u, \vct u \rangle_{\symbfcal{H}} \geq \gamma \langle \vct u, \vct u \rangle_{\symbfcal{H}}.
  \end{equation}
\item A variant of the Sobolev embedding theorem yields the embedding \(H^1_{\nu}(\R, \symbfcal{H}) \hookrightarrow C_{\nu}(\R, \symbfcal{H})\) (cf.~\autocite[Lemma 3.1.59]{picardPartialDifferentialEquations2011}), where
  \[
    C_{\nu}(\R, \symbfcal{H}) \coloneq \{ \vct f \colon \R \to \symbfcal H \mid \vct f \text{ is continuous}, \, \sup_{t \in \R} |\vct f(t)| \e^{-\nu t} < \infty \}\,.
  \]
\item If \(\vct F \in H_\nu^k(\R; \symbfcal{H})\) for \(k \in \N\), then \(\vct U \in H_\nu^k(\R; \symbfcal{H})\) and~\eqref{eq:evolutionary-problem} is solved literally, i.e., \((\partial_t \symbfcal M_0 + \symbfcal M_1 + \symbfcal A) \vct U = \vct F\). Moreover, \(\vct U \in C_{\nu}(\R, \dom(\symbfcal{A}))\).
\item To address initial value problems, the solution theory must be extended to
  allow $\vct u$ to satisfy initial conditions. This approach is discussed in
  detail in~\cite{trostorffExponentialStabilityInitial2018b} and~\autocite[Section
  6.3]{seifertEvolutionaryEquationsPicard2022}.
\end{itemize}
\subsection{Well-posedness of nonlinear evolutionary
  problems\label{ssec:wp-nlevolutionary-problem}}
In this section, we review nonlinear perturbartions of Maxwell's equations~\eqref{eq:Maxwell_pos_init} as
done in~\cite{dohnalWellPosednessExponentialStability2023}, where the authors
employ a Banach fixed-point methodology to equations of the form
\begin{equation}\label{eq:nl-evolutionary-problem}
(\partial_t\symbfcal M(\partial_t) + \symbfcal A)\vct u = \vct f(\vct u)\,,
\end{equation}
by formulating~\eqref{eq:nl-evolutionary-problem} as a fixed-point equation
$\vct u=\symbfcal{S}_{\nu}\vct f (\vct u)$.
\begin{definition}{Uniformly Lipschitz-continuous map}{uniform-lip-map}
Let \( \vct f \) be a function defined on \( \dom(\vct f) \subseteq \bigcap_{\nu \ge \nu_0}L^2_\nu(\R,\,\symbfcal{H}) \). The function \( \vct f \) is called \emph{uniformly Lipschitz-continuous} if the following conditions are met:
\begin{itemize}\itemsep1pt \parskip0pt \parsep0pt
    \item The domain \( \dom(\vct f) \) is dense in \( L^2_\nu(\R,\,\symbfcal{H}) \).
    \item The function \( \vct f \) can be extended to a Lipschitz-continuous map, denoted as \( \vct f_\nu \), mapping \( L^2_\nu(\R,\,\symbfcal{H}) \) to itself, such that when \( \nu \) and \( \nu' \) are both greater than \( \nu_0 \), the functions \( \vct f_\nu \) and \( \vct f_{\nu'} \) coincide on the intersection \( L^2_\nu(\R,\,\symbfcal{H})\cap L^2_{\nu'}(\R,\,\symbfcal{H}) \).
\end{itemize}
For any \( \nu > \nu_0 \), it is permissible to represent this function simply as \( \vct f\colon L^2_{\nu}(\R,\,\symbfcal{H}) \to L^2_{\nu}(\R,\,\symbfcal{H}) \).
\end{definition}
\paragraph*{Assumptions required for the fixed point arguments}
The additional assumptions to Theorem~\ref{thm:picards} for the well-posedness of~\eqref{eq:nl-evolutionary-problem} are as follows:
\begin{description}\itemsep1pt \parskip0pt \parsep0pt
\item[\customlabel{itm:lipschitz-perturb1}{(A1)}] For some constant $d>0$, the material law \(\symbfcal{M}(\partial_t)\) satisfies the condition
    \begin{equation}\label{eq:Lipschitz_perturb}
    \re(z\symbfcal{M}(z)) \geq \frac{\re z}{d}, \quad \text{for all } z\in \C
    \text{ with }\re z > \nu_0.
    \end{equation}
  \item[\customlabel{itm:lipschitz-perturb2}{(A2)}] The nonlinearity \(\vct{f}\) is \emph{uniformly Lipschitz continuous}
      for all \(\nu > \nu_0\), i.\,e.\
    \begin{equation}\label{eq:Lipschitz_perturb-i}
    \|\vct{f}(\vct{u}) - \vct{f}(\vct{v})\| \leq L\|\vct{u} - \vct{v}\| \quad \text{in } L^2_\nu(\mathbb{R},\,\symbfcal{H}),
    \end{equation}
    and \(\limsup_{\nu \to +\infty} \frac{d}{\nu} \|\vct{f}\|_{\mathrm{Lip}}=\sup_{\vct{u} \neq \vct{v} \in \dom(\vct f)} \frac{\|\vct{f}(\vct{u}) - \vct{f}(\vct{v})\|}{\|\vct{u} - \vct{v}\|} < 1\).
\end{description}
Now we give an abstract result for the well-posedness of a nonlinear first
order evolutionary system.
\begin{lemma}{Well-posedness of nonlinearly perturbed equations~\autocite[Proposition 2.8]{dohnalWellPosednessExponentialStability2023}}{wp-nl-perturb}
Given the assumptions~\ref{itm:lipschitz-perturb1}
and~\ref{itm:lipschitz-perturb2}, there exists a \(\nu_1 \geq \nu_0\) such
that for every \(\nu > \nu_1\), the nonlinear evolutionary equation
\begin{equation}
  \label{eq:wp-nlpe}
  (\partial_t \symbfcal{M}(\partial_t) + \symbfcal{A}) \vct{u} = \vct{f}(\vct{u})
\end{equation}
    admits a unique solution in \(L^2_\nu(\mathbb{R}, \symbfcal{H})\) that does not depend on \( \nu \).
  \end{lemma}
  Lemma~\ref{lem:wp-nl-perturb} is applied to the Maxwell equation~\eqref{eq:Maxwell_pos_init} with
  nonlinear polarization.
\paragraph*{Nonlinear Maxwell system}
The nonlinearly perturbed equation~\eqref{eq:wp-nlpe} is adapted in accordance
with Maxwell equation~\eqref{eq:Maxwell_pos_init} by setting the electric
displacement vector to
\[
\vct{D}(\vct{E}) = \varepsilon(\partial_t) + \vct{P}_{\mathrm{nl}},
\]
where $\varepsilon(\partial_t)$ represents the linear material law and \(\vct{P}_{\mathrm{nl}}\) is the nonlinear polarization term
\begin{equation}\label{eq:pol-nl}
    \vct{P}_{\mathrm{nl}}(\vct{E}) = \int_{\mathbb{R}} \mat K(t - s)
    \vct{q}(\vct{E}(s)) \drv s\,,
\end{equation}
with $\mat K\colon \R\to \symbfcal{B}(\symbfcal{H})$, $\supp \mat K
\subseteq [0,\,\infty)$ and Lipschitz continuous $\vct q \colon \symbfcal{H}\to \symbfcal{H}$.
The \emph{nonlinear Maxwell system} reads
\begin{equation}\label{eq:Maxwell_nonlinear}
    \left( \partial_t \begin{pmatrix} \varepsilon(\partial_t) & 0 \\ 0 & \mu_0 \end{pmatrix} + \begin{pmatrix} 0 & -\curl \\ \curl_0 & 0 \end{pmatrix} \right)
    \begin{pmatrix} \vct{E} \\ \vct{H} \end{pmatrix} = \begin{pmatrix} \vct f -\partial_t \vct{P}_{\mathrm{nl}}(\vct{E}) \\ 0 \end{pmatrix}\,.
  \end{equation}
  The uniform Lipschitz continuity of~\eqref{eq:pol-nl} is needed
  for Lemma~\ref{lem:wp-nl-perturb} to be applicable to~\eqref{eq:Maxwell_nonlinear}.
  This relies on the continuity and differentiability of
$\mat K$ and Lipschitz continuity of $\vct q$~\autocite[Lemma 2.10]{dohnalWellPosednessExponentialStability2023}.
\paragraph*{Well-posedness of the nonlinear Maxwell system with saturable nonlinearity}
To include the most prevalent models in nonlinear optics, the conditions
needed for Lemma~\ref{lem:wp-nl-perturb} are too restrictive. The authors
of~\cite{dohnalWellPosednessExponentialStability2023} propose \emph{saturable
nonlinearities} of the form
\[\vct q(\vct u)(\vct x) = \frac{|\vct u(\vct x)|^{k-1}}{1 + \tau |\vct u(\vct
    x)|^{k-1}} \vct u(\vct x) \coloneq V(|\vct u(\vct x)|) \vct u(\vct x)\,,\]
where $k\geq 2$, $\tau >0$, in order to bridge the gap between modeling and
mathematical analysis. For $k=2$ this represents a saturable second order
nonlinearity and for $k=3$ this represents a saturable modification of the
Kerr-type nonlinearity. Since the mapping $\R^3\ni \xi \mapsto V(\abs{\xi})\xi$
is smooth and asymytotically linear, it is Lipschitz continuous, which $\vct q$
inherits. It therefore fulfills the requirement of
Lemma~\ref{lem:wp-nl-perturb}. Thus, the~\emph{nonlinear Maxwell
  system}~\eqref{eq:Maxwell_nonlinear} with saturable nonlinearity is
well-posed.

\subsubsection*{Local well-posedness of the nonlinear Maxwell equation}
\label{sec:org2babe62}
The requisite of uniform Lipschitz-continuity within the range of spaces
\(L^2_{\nu}(\mathbb{R},\,\symbfcal{H})\) is constraining, especially since it
excludes the consideration of nonlinearities that exhibit superlinear growth as
potential candidates for \(\vct q\). This requirement can be relaxed by
substituting it with Lipschitz continuity on closed subsets within
\(L^2_{\nu}(\mathbb{R},\,\symbfcal{H})\). For sufficiently large values of
\(\nu\), these subsets then include the given data~\autocite[Proposition
2.12]{dohnalWellPosednessExponentialStability2023}. We review this more nuanced
version of the fixed-point approach for establishing local well-posedness.

\paragraph*{Assumptions required for the refined fixed point arguments.}
We begin by laying down the necessary assumptions for local well-posedness:
\begin{description}\itemsep1pt \parskip0pt \parsep0pt
\item[\customlabel{itm:lipschitz-perturb1p}{(A1\('\))}] The material law $\symbfcal{M}(\partial_t)$ satisfies~\eqref{eq:Lipschitz_perturb}.
\item[\customlabel{itm:lipschitz-perturb2p}{(A2\('\))}] The nonlinearity \(\vct{f}\colon L^2_{\nu}(\mathbb{R}, \symbfcal{H})
      \to L^2_{\nu}(\mathbb{R}, \symbfcal{H})\) is locally Lipschitz continuous
      on a closed subset \(\mathcal{D} \subset L^2_{\nu}(\mathbb{R},
      \symbfcal{H})\), meaning for some constants \(c > 0\) and \(\alpha > 0\):
      \begin{equation}
        \label{eq:lipschitz-refined-i}
\lVert \vct f(\vct u) - \vct f(\vct v) \rVert_{L^2_{\nu}(\mathbb{R}, \symbfcal{H})} \le c \left( \lVert \vct u \rVert_{L^2_{\nu}(\mathbb{R}, \symbfcal{H})} + \lVert \vct v \rVert_{L^2_{\nu}(\mathbb{R}, \symbfcal{H})} \right)^\alpha \lVert \vct u - \vct v \rVert_{L^2_{\nu}(\mathbb{R}, \symbfcal{H})}, \quad \forall \vct u, \vct v \in \mathcal{D}.
      \end{equation}
    \item[\customlabel{itm:lipschitz-perturb3p}{(A3\('\))}] The forcing term \(g \in L^2_{L^2_{\nu}(\mathbb{R}, \symbfcal{H})}(\mathbb{R}, \symbfcal{H})\) satisfies \(\lVert g \rVert_{L^2_{\nu}(\mathbb{R}, \symbfcal{H})} = o\left( \nu^{1 + \frac{1}{\alpha}} \right)\) as \(\nu \to \infty\).
\end{description}

\begin{proposition}{Fixed point refinement~\cite{dohnalWellPosednessExponentialStability2023}}{fpr}
Let \(\symbfcal A: \dom(\symbfcal A) \subset \symbfcal{H} \to \symbfcal{H}\) be a skew self-adjoint operator, and \(\symbfcal M\) denote a linear material law satisfying~\ref{itm:lipschitz-perturb1p}.
Let \(\vct f: L^2_{\nu}(\mathbb{R},\,\symbfcal{H}) \to L^2_{\nu}(\mathbb{R},\,\symbfcal{H})\) be a causal nonlinear map with \(\vct f(0) = \vct 0\), satisfying~\ref{itm:lipschitz-perturb2p}
and assume \(\vct g \in L^2_{\nu}(\mathbb{R},\,\symbfcal{H})\) satisfies~\ref{itm:lipschitz-perturb3p}.
Then, for \(\nu > \nu_0\), the nonlinear system:
\[
  (\partial_t\symbfcal M(\partial_t) + \symbfcal A)\vct u = \vct f(\vct u) + \vct g
\]
admits a unique solution \(\vct u \in L^2_{\nu}(\mathbb{R},\,\symbfcal{H})\).
\end{proposition}

\begin{proofp}
Define the linear solution operator \(\symbfcal S_\nu\) as follows:
\[
\symbfcal S_\nu \coloneq \left( \overline{\partial_t\symbfcal M(\partial_t) + \symbfcal A} \right)^{-1}.
\]
Given \(\lVert \symbfcal S_\nu \rVert \le \frac{d}{\nu}\), we apply a fixed-point argument to show that \(\symbfcal S_\nu (\vct f(\vct u) + \vct g)\) is a contraction on a closed ball \(B_r \subset L^2_{\nu}(\mathbb{R}, \symbfcal{H})\). For sufficiently large \(\nu\), the Lipschitz constant of \(\vct{f}\) ensures the contraction condition is satisfied, guaranteeing a unique solution \(\vct u\).
\end{proofp}

\paragraph*{Fixed-point arguments for nonlinearities with superlinear growth}
To apply the fixed point methodology to a broader class of nonlinearities, we
adapt the fixed-point argument to nonlinearities that are subject to temporal
cutoffs.
\begin{definition}{Nonlinearity with temporal cutoff}{nltc}
Let \(\vct f_T\) be a truncated version of \(\vct f\) with a temporal cutoff:
\[
\vct f_T(\vct u)(t) = \int_{\mathbb{R}} \int_{\mathbb{R}} \mat K_T(t, t - \tau_1, t - \tau_2) \vct q(\vct u(\tau_1), \vct u(\tau_2)) \, \drv \tau_1 \, \drv \tau_2,
\]
where \(\mat K_T\) is a kernel truncated at \(T > 0\). This modification allows for localized control over nonlinear growth within the time interval \([0, T]\).
\end{definition}

By introducing a cutoff in the kernel \(\mat{K}_T\), we ensure that the nonlinearity remains Lipschitz continuous within the truncated time interval. This enables us to apply the fixed-point method even for nonlinearities with superlinear growth.

\subsubsection*{Local well-posedness for nonlocal quadratic polarization}

To illustrate the local well-posedness of the nonlinear Maxwell equation, we
consider the Maxwell system with a second-order nonlocal quadratic polarization.

\begin{definition}{Fully nonlocal quadratic polarization}{fnlqp}
A fully nonlocal quadratic polarization \(\vct{p}^{(2)}(\vct{e})\) is defined
as:
\begin{equation}
  \label{eq:fnlqp}
  \vct{p}^{(2)}(\vct{e})(t) = \int_{\mathbb{R}} \int_{\mathbb{R}} \mat K(t - \tau_1, t - \tau_2) \vct{q}(\vct{e}(\tau_1), \vct{e}(\tau_2)) \, \drv \tau_1 \, \drv \tau_2,
\end{equation}
where \(\mat K\) is a causal operator-valued kernel with compact support, and \(\vct{q}\) is a bounded bilinear map.
\end{definition}

\paragraph*{Assumptions on the integrability of $\mat K$}
To fulfil the assumption~\ref{itm:lipschitz-perturb2p}, we assume that the kernel \(\mat K\) satisfies the integrability conditions
\begin{equation}
  \label{eq:integrability-k}
  \begin{aligned}
 L^2_{\mat K} &= \iint_{0}^{\infty} \lVert \mat K(\tau_1, \tau_2) \rVert \e^{-\nu_{\mat K}(\tau_1 + \tau_2)} \, \drv\tau_1 \, \drv\tau_2 < \infty\,, \\
  \ell_{\mat K} &= \sup_{\tau_1, \tau_2 \in \R} \int_{0}^{\infty} \lVert \mat K(t-\tau_1, t-\tau_2) \rVert \e^{-\nu_{\mat K}(2t-\tau_1-\tau_2)} \, \drv t < \infty,
  \end{aligned}
\end{equation}
for some \(\nu_{\mat K} \in \R\).

As a consequence of the assumptions~\eqref{eq:integrability-k}, using the bilinearity of \(\vct{q}\), we have the inequality
\[
  \lVert \vct{f}_T(\vct{u}) - \vct{f}_T(\vct{v}) \rVert_{L^2_{\nu}(\mathbb{R}, \symbfcal{H})} \le
  \sqrt{T}\e^{\nu T}C_{\vct{q}} \sqrt{d_{\mat K}L_{\mat K}} \left( \lVert
    \vct{u} \rVert_{L^2_{\nu}(\mathbb{R}, \symbfcal{H})} + \lVert \vct{v} \rVert_{L^2_{\nu}(\mathbb{R}, \symbfcal{H})} \right) \lVert \vct{u} -
  \vct{v} \rVert_{L^2_{\nu}(\mathbb{R}, \symbfcal{H})},
\]
where \(\vct{f}_T\) represents the nonlinearity with a temporal cutoff.
The time derivative of $\vct{p}^{(2)}$ preserves the nonlocal quadratic nonlinearity and can be computed as
\[
\partial_t\vct{p}^{(2)}(\vct{e})(t) = \int_{\mathbb{R}} \int_{\mathbb{R}} (\partial_1 + \partial_2)\mat K(t - \tau_1, t - \tau_2) \vct{q}(\vct{e}(\tau_1), \vct{e}(\tau_2)) \, \drv \tau_1 \, \drv \tau_2\,.
\]
Here, the assumptions~\eqref{eq:integrability-k} ensure the boundedness and integrability of $\partial_t\vct{p}^{(2)}$.
\begin{proposition}{Local well-posedness of Maxwell system with quadratic polarization}{lwpmwqp}
  For sufficiently large \(\nu\), the Maxwell system with fully nonlocal quadratic polarization~\eqref{eq:fnlqp} admits a unique solution \(\vct u \in L^2_{\nu}(\mathbb{R},\,\symbfcal{H})\), provided that the initial data is sufficiently small.
\end{proposition}

\begin{proofp}
The proof follows from the fixed-point refinement and the temporal cutoff for the nonlinearity, ensuring that the nonlocal quadratic polarization remains well-posed within the given time interval.
\end{proofp}
\paragraph*{Instantaneous nonlinearities}
So far instantaneous (i.\,e.\ zero-delay) nonlinearities are excluded from the
right-hand side of the system~\eqref{eq:nl-evolutionary-problem} by imposing the
condition \(\supp \chi^{(2)} \subset (0, \infty)^2\). This condition can be
relaxed by working in the Sobolev space \(H^1_{-\nu}\) and utilizing the Sobolev
embedding
$H^1_{\nu}(\R,\,\symbfcal{H})\hookrightarrow C_{\nu}(\R,\,\symbfcal{H})$, to
derive the integrability conditions for the perturbation arguments. However, we
note that higher derivatives pose additional challenges. For instance,
estimating terms such as
\[
\int_{\mathbb{R}} \chi^{(2)}(0, \tau) \partial_t \vct Q(\vct{u}(t), \vct{u}(t - \tau)) \drv\tau
\]
in the \(H^1_\rho\)-norm requires more regularity of \(\vct{u}\) to handle the
higher-order derivative terms that arise in the analysis. Therefore, while
working in \(H^1_{-\nu}\) mitigates the need to exclude instantaneous
nonlinearities, it introduces a requirement for higher regularity in the
solutions.

Removing spatial dispersion from the nonlinearity leads to similar problems. If
\( \vct Q\) is a matrix-valued bilinear operator, working in \( H^k\)-Sobolev
spaces with $k>\frac{d}{2}$ is desired due to their algebra property. However,
additional spatial regularity cannot typically be inferred from regular data.
These remarks show that simpler, instantaneous, and local nonlinear material
laws introduce additional problems that usually require more regularity of the
solution. Overall, quasilinear systems are difficult to handle in the
evolutionary framework due to the low regularity.

In summary, the local well-posedness of nonlinear Maxwell equations is achieved
through a refined fixed-point approach and temporal cutoff, making it applicable
to a broader class of nonlinearities, including those in nonlinear optics. These
theoretical foundations will support the development of computational
electromagnetics models in future work. However, due to the low regularity of
the Hilbert space approach, unavoidable restrictions must be acknowledged.

\subsection{Exponential stability\label{sec:exp-stab}}
\label{sec:org6980be3}
Exponential stability means that the solution of an initial value problem decays
at an exponential rate as time approaches infinity. In this context,
exponentially decaying right-hand sides result in exponentially decaying
solutions. Due to the lack of continuity with respect to time, a pointwise
definition is impractical. Instead, we use exponentially weighted spaces
\(L^2_{\nu}(\R; \symbfcal{H})\) with negative \(\nu\) to define exponential
stability by requiring the invariance of these spaces under the solution
operator. Exponential stability is crucial for dispersion models and for
ensuring the effectiveness of PML\@. We use the definition
in~\autocites{trostorffExponentialStabilityLinear2013a}[Section
11.1]{seifertEvolutionaryEquationsPicard2022}.
\begin{definition}{Exponential stability}{exp-stab}
  Let \(\symbfcal{H}\) be a Hilbert space, $\symbfcal A$ be a skew self-adjoint
  operator, and let $\symbfcal M$ represent a linear material law. Suppose
  $z\symbfcal M(z)$ is strictly accretive on a half-plane
  $\C_{\re > \nu_0}$ for some $\nu_0 \in \R$ such that~\eqref{eq:MaterialLawAccretivity} holds. The problem is \emph{exponentially stable with decay rate
    \(\rho_0 > 0\)} if, for all \(\rho > [0,\,\rho_0)\), $\nu \geq \nu_0$ and
  \(\vct f \in L_{-\rho}^2(\mathbb{R}; \symbfcal{H}) \cap L_{\nu}^2(\mathbb{R};
  \symbfcal{H})\),
  \begin{equation*}
    \left(\overline{\partial_{t,\nu} \symbfcal{M}(\partial_{t,\nu}) + \symbfcal{A}}\right)^{-1} \vct f \in L_{-\rho}^{2}(\mathbb{R}; \symbfcal{H}).
  \end{equation*}
\end{definition}
For homogeneous materials, exponential stability is achieved if the electric
permittivity \(\varepsilon(\partial_t)\) and \(\mu\) satisfy specific
conditions.
\begin{remark}{Conditions for exponential
    stability~\autocite{dohnalWellPosednessExponentialStability2023}}{maccret}
  \begin{description}
  \item[\customlabel{itm:MatCond_accr1}{(M2)}] For all \(\delta > 0\), there exist
    \(\nu > 0\) and \(c > 0\) such that
    \(\forall z \in \C_{\re > -\nu} \setminus B[0,\delta]: \; \re z
    \varepsilon(z) \ge c\).
  \item[\customlabel{itm:MatCond_accr2}{(M3)}]
    \(\varepsilon(\partial_t) = \varepsilon_0 + \chi(\partial_t)\), where
    \(\varepsilon_0\) is bounded, linear, self-adjoint, uniformly positive
    definite, \(\lim_{z \to 0} z \chi(z) = 0\), and \(\chi(z)\) and
    \(z \chi(z)\) are bounded in \(\C_{\re > -\nu_1}\), with
    \(\varepsilon_0 + \re \chi(z) \ge c_1 > 0\).
  \end{description}
\end{remark}

Under these assumptions, \(\varepsilon(z) = \varepsilon_0 + z^{-1}(z \chi(z))\)
satisfies the conditions of Definition~\ref{def:exp-stab} with
\(\symbfcal{M}_0 = \varepsilon_0\) and
\(\symbfcal{M}_1(z) = z \chi(z)\)~\autocite[Theorem
3.10]{dohnalWellPosednessExponentialStability2023}. The
conditions~\ref{itm:MatCond_accr1} and~\ref{itm:MatCond_accr2} are
sufficient for exponential stability for the second-order system
\(\symbf\nabla \times \symbf\nabla \times \vct{E} + \mu_0 \partial_{tt} \vct{D}
+ \partial_{t} \vct{J} = 0\) on a bounded domain \(\Omega\)~\autocite[Theorem
3.15]{dohnalWellPosednessExponentialStability2023}.
Alternatively, stronger material damping can be imposed for exponential
stability, as in~\autocite[Proposition
2.1.5]{trostorffExponentialStabilityInitial2018b}.

\begin{proposition}{A Condition for exponential stability}{exp-stab-pos-def}
  Let \(\symbfcal{H}\) be a Hilbert space, \(\symbfcal{A}\) a skew self-adjoint
  operator, and \(\symbfcal{M}_0, \symbfcal{M}_1\) bounded, analytic, and linear
  operators. Assume there exists \(\nu_0 > 0\) such that
  \(\C_{\re > -\nu_0} \setminus \dom(\symbfcal{M}_0)\) is discrete and
  \begin{equation}
    \label{exp-stab-cond}
    \exists c > 0 \;\forall z \in \dom(\symbfcal{M}) \cap \C_{\re > -\nu_0},\; u \in \symbfcal{H} : \re \dd{z \symbfcal{M}(z) \vct{u}}{\vct{u}}_{\symbfcal{H}} \geq c \norm{\vct{u}}_{\symbfcal{H}}^{2}.
  \end{equation}
  Then the problem associated with \(\symbfcal{M}\) and \(\symbfcal{A}\) is
  well-posed and exponentially stable with decay rate \(\nu_0\).
\end{proposition}

For conductive materials, we use another set of conditions for
exponential stability~\autocite[Section
4.1]{dohnalWellPosednessExponentialStability2023}.

\begin{remark}{Conditions for exponential
    stability~\autocite{dohnalWellPosednessExponentialStability2023}}{maccret2}
  \begin{description}
  \item[\customlabel{itm:MatCond_accr_conductivity}{(M2\('\))}]
    \(\re z \varepsilon(z) \ge 0\) for \(\re z > 0\).
  \item[\customlabel{itm:MatCond_conductivity}{(M3\('\))}]
    \(\varepsilon(\partial_t) = \varepsilon_0 + \chi(\partial_t)\) with
    \(\varepsilon_0\) self-adjoint and uniformly positive definite,
    \(\lim_{z \to 0} z \chi(z) = 0\), and \(\chi(z)\) and \(z \chi(z)\) bounded
    on \(\C_{\re > -\nu_1}\), with \(\varepsilon_0 + \re \chi(z) \ge c_1 > 0\).
  \end{description}
\end{remark}

Under these conditions, the material law
\(\symbfcal{M}(z) = \varepsilon(z) + z^{-1} \sigma\) satisfies for some
\(\nu_0 > 0\), \(c > 0\),
\begin{equation}
  \re z \ge -\nu_0 \implies \re z \symbfcal{M}(z) = \re z \varepsilon(z) + \sigma \ge c.
  \label{eq:accr_perm_strict}
\end{equation}

\section{A class of dispersive Maxwell equations\label{sec:dispersion-picard}}
In this study, we have analyzed classes of nonlinear models based on their
mathematical properties and demonstrated that these results apply to practically
relevant physical systems. This analysis also covers general linear material
laws (cf. Definition~\ref{def:matlaw}). However, we did not address specific
physical models so far. In order to establish well-posedness and exponential
stability for physically relevant models, we now consider a generalized
dispersion model as introduced
in~\cite{viqueratSimulationElectromagneticWaves2015}. This model encapsulates a
broad spectrum of prevalent dispersion models, notably including the Debye and
Lorentz models. The Lorentz model describes material resonance in response to
electromagnetic fields and is particularly relevant in high-frequency domains,
such as optics. In contrast, the Debye model captures the low-frequency
dielectric relaxation process, whereby dipoles in a material reorient in
response to a changing electric field. Its applications are most significant in
lower frequencies, such as in microwave or radio wave regimes, where the
time-dependent alignment of molecules with the field is of significance.
\begin{definition}{Generalized dispersion model}{gen-disp-model}
Let  $L_1$ and $L_2$ be disjoint index sets and $\varepsilon_{\omega}$, $\sigma$,
$(a_{l})_{l\in L_1}$, $(b_{l})_{l\in L_1}$, $(c_{l})_{l\in L_2}$,
$(d_{l})_{l\in L_2}$, $(e_{l})_{l\in L_2}$, $(f_{l})_{l\in L_2}$ real
constants. Then we define the generalized dispersion models in the
Fourier-Laplace domain by
\begin{equation}
  \label{eq:disp-mod}
  \varepsilon(z) =  \varepsilon_{\omega}
  + \frac{\bar{\sigma}}{z}
  + \sum_{i\in L_{1}} \frac{a_{l}}{b_{l}+ z}
  + \sum_{l\in L_{2}} \frac{c_{l}+ z d_{l}}{e_{l}+z f_{l}+z^{2}}\,.
\end{equation}
The constant $\varepsilon_{\omega}$ is the refractive index in the high frequency limit and
$\sigma$ the conductivity.
\end{definition}
By examining the equation~\eqref{eq:disp-mod}, it is evident that it includes
the Debye model via the first-order terms over the index set $L_1$. The Lorentz
model is represented by the terms associated with $L_2$. In the time domain, we
denote the generalized dispersion model~\eqref{eq:disp-mod} by
$\varepsilon(\partial_t)$, according to~\ref{def:matlaw}.

We formulate auxiliary differential equations (ADEs) in order to incorporate the
generalized dispersion model~\eqref{eq:disp-mod} into the Maxwell
equations~\eqref{eq:Maxwell_pos_init}. For a detailed description of the
derivation of these ADEs we refer
to~\cite{margenbergAccurateSimulationTHz2023,margenbergOptimalDirichletBoundary2024}.
We then assess the well-posedness of the resulting set of equations, according
to Theorem~\ref{thm:picards}. This approach leads to the derivation of the
following system of equations
\begin{subequations}
  \label{eq:disp-pde}
  \begin{align}
    \varepsilon_{\omega} \partial_{t}\vct E + \paran{\bar{\sigma} + \sum_{l\in L_{2}}d_{l}}\vct E +
    \sum_{l\in L_{1}} a_l\vct E-b_l\vct p_l +
    \sum_{l\in L_{2}} \vct j_{l}-\curl \vct H&=\vct f\\
    \partial_{t}\vct H + \curl_0 \vct E &= 0\\
    \partial_t \vct p_{l} - a_{l} \vct E + b_{l} \vct p_{l} &= 0\,,\quad l\in L_{1}\\
    \partial_{t} \vct j_{m} + (f_md_m-c_m)\vct E + f_{m}\vct j_{m} + e_{m}\vct
    p_{m}=0\,,\quad \partial_{t} \vct p_{m} - d_{m} \vct E - \vct j_{m} &=
                                                                           0\,,\quad
                                                                           m\in L_{2}\,.
  \end{align}
\end{subequations}
These equations can be cast into the form required by Picard's
theorem. To simplify notation we introduce the notation of
\(\vec(\bullet)\) and \(\diag(\bullet)\) which expands a given set into a column
vector or diagonal matrix respectively. We further define
\(x_{\symcal{S}}=\{x_{s}\;:\; s\in \symcal{S}\}\) to denote a sequence over an
index set \(\symcal{S}\). We define \(\rho\coloneq\sum_{l\in L_2} d_l+\sum_{l\in
  L_1 a_l}\) and put
\begin{subequations}\label{eq:disp-ops}
\begin{equation}\footnotesize
  \symbfcal A=\begin{tikzpicture}[baseline=(current bounding box.center)]%
      \matrix[sparsemat] {
         &-\curl &\cdots & \\
        \curl_0& &\cdots  &  \\
        \vdots&\vdots   & \ddots & \vdots \\
        &   &\cdots &   \\
      };
    \end{tikzpicture}\;,\quad
    \symbfcal{M}_0=\begin{tikzpicture}[baseline=(current bounding box.center)]%
      \matrix[sparsemat] {
        \varepsilon_{\omega} & &\cdots & \\
        & 1 &\cdots  &  \\
        \vdots&\vdots   & \ddots & \vdots \\
        &   &\cdots &   1\\
      };
    \end{tikzpicture}
    \normalsize
    \end{equation}\;,
    \begin{equation}\footnotesize
    \symbfcal{M}_1=\begin{tikzpicture}[baseline=(current bounding
             box.center)]%
             \matrix[sparsemat] {
               \bar{\sigma} + \rho &         & -\vec(b_{L_1})^{\top}  & 1&  \\
                            & &         &              & \\
               -\vec(a_{L_{1}})  & & \diag(b_{L_1})       &  &     \\
               \vec(d_{L_2}f_{L_2})-\vec(c_{L_2}) & &         & \diag(f_{L_2}) & \diag(e_{L_2})\\
               -\vec(d_{L_2})  & & & -1 & \\
             };
           \end{tikzpicture}\:.
    \normalsize
\end{equation}
\end{subequations}
and
\begin{equation}
  \label{eq:disp-sols}
  \vct U={(\vct E,\,\vct H,\,\vct p_{L_{1}},\,\vct j_{L_{2}},\,\vct p_{L_{2}})}^{\top},\; \vct f={(\vct f,\,0,\,0,\,0,\,0)}^{\top}.
\end{equation}
\begin{problem}{Dispersive Maxwell equation as an evolutionary
    problem}{evolution-disp}
  Let \(\symbfcal{H} = L^d \times L^d\), a Hilbert space with the inner product
  of \(L^{6d}\). Let
  $\symbfcal M_0,\: \symbfcal M_1\colon \symbfcal{H}\to \symbfcal{H}$ be defined
  by~\eqref{eq:disp-ops}. The evolutionary problem with $\vct f$ and $\vct U$
  defined according to~\eqref{eq:disp-sols} reads as:
  \emph{For a given right-hand side $\vct f\in L^2_{\nu}(\R;\:\symbfcal{H})$ find
    $\vct U\in L^2_{\nu}(\R;\:\symbfcal{H})$ such that}
  \begin{equation}
    \label{eq:evo-problem-disp}
    (\partial_t \symbfcal M_0 + \symbfcal M_1+\symbfcal A)\vct U=\vct f.
  \end{equation}
\end{problem}
\begin{corollary}{Well-posedness of
    Problem~\ref{problem:evolution-disp}}{well-posed-evo-disp}
  Problem~\ref{problem:evolution-disp} is well-posed. That is for each
  $\vct f\in L^2_{\nu}(\R;\:\symbfcal{H})$ there exists a unique solution $\vct U\in L^2_{\nu}(\R;\:\symbfcal{H})$.
\end{corollary}
\begin{proofp}
  For Problem~\ref{problem:evolution-disp} operator $\symbfcal A$ is skew self-adjoint. The
  operators $\symbfcal M_0$, $\symbfcal M_1$ are clearly bounded linear operators and $\symbfcal M_0$ is
  self-adjoint and positive definite.
  We now check, whether there exist $\gamma >0$ and
  $\nu_{0}>0$ such that condition~\eqref{eq:evo-cond}
  is fulfilled.
  Since $\symbfcal M_0$ is positive definite, we have
\begin{align*}
    \label{wp-ineq}
    \re \dd{z \symbfcal M(z) \vct u}{\vct u}_{\symbfcal{H}} &= \re \dd{z \symbfcal M_0 \vct u}{\vct u}_{\symbfcal{H}} + \dd{\symbfcal M_1 \vct u}{\vct u}_{\symbfcal{H}} \notag \\
    &> -\nu \min\{\varepsilon_{\omega},\,1\}\dd{\vct u}{\vct u}_{\symbfcal{H}} + \dd{\symbfcal M_1 \vct u}{\vct u}_{\symbfcal{H}} \,.
  \end{align*}
  As $\symbfcal M_1$ is linear and bounded, there exists a $\nu_0>0$ such that
  negative terms arising from $\re \symbfcal M_{1}$ can be absorbed and the
  inequality in~\eqref{eq:evo-cond} with a $\gamma >0$ holds. A (not necessarily
  optimal) lower bound for $\nu_0$ can be found in a straightforward manner by
  applying Young's inequality to $\dd{\symbfcal{M}_1\vct u}{\vct u}$. Then, we
  can find a $\gamma>0$ such that~\eqref{eq:evo-cond} holds. Therefore,
  Theorem~\ref{thm:picards} concludes the proof.
\end{proofp}

\subsection*{Exponential stability of different dispersion models}
First, we
show exponential stability of the Lorentz model and the Debye model. By
linearity, exponential stability of these models extend to the full generalized
dispersion model. Upon initial examination, the material law represented by
\(\symbfcal M(\partial_t) = \Bigl(\begin{smallmatrix} \varepsilon(\partial_t) &
  0 \\ 0 & \mu \end{smallmatrix}\Bigr)\) does not meet any of the strict
accretivity conditions~\ref{itm:MatCond_accr1}
and~\ref{itm:MatCond_accr2} (cf.~Remark~\ref{rem:maccret}). This is evident
from
\begin{equation}\label{eq:no-accretivity}
\re \dds{z \symbfcal M(z) \begin{pmatrix}{\color{gray} \symbf 0} \\ \vct u \end{pmatrix}}{ \begin{pmatrix}{\color{gray} \symbf 0} \\ \vct u \end{pmatrix}} = (\re z)\dd{\mu \vct u}{\vct u} \rangle = 0\,,
\end{equation}
which holds true whenever \(\re z=0\), irrespective of the value of
\(\varepsilon\). While we cannot conclusively derive exponential stability, a
potential avenue we consider here is to use the results
in~\autocite{dohnalWellPosednessExponentialStability2023}, in particular Theorem
3.10 and 3.15 in conjunction with~\ref{itm:MatCond_accr1} and
\ref{itm:MatCond_accr2} of Remark~\ref{rem:maccret} to obtain the
exponential stability. For the conductivity $\bar \sigma$, it is straightforward
to show exponential stability, see also the proof of
Corollary~\ref{cor:well-posed-evo-b}.
\begin{corollary}{Exponential stability of the Debye model}{exp-stab-deb}
Consider Problem~\ref{problem:evolution-disp} and suppose $\bar{\sigma}=0$,
$L_{1}=\brc{0}$, $L_{2}=\emptyset$. Equation~\eqref{eq:disp-mod}
is equivalent to a Debye model, i.\,e.\
\begin{equation}
  \label{eq:abstract-material-deb}
  \varepsilon(z) =  \varepsilon_{\omega}+\underbrace{\frac{a}{z + b}}_{\coloneq\chi(z)},\quad z\in \C\,.
\end{equation}
The associated material law is exponentially stable as given by Definition~\ref{def:exp-stab}.
\end{corollary}
\begin{proofp}
  We show the exponential stability of~\eqref{eq:abstract-material-deb} by using
  the criteria~\ref{itm:MatCond_accr_conductivity}
  and~\ref{itm:MatCond_conductivity} (cf.~Remark~\ref{rem:maccret2}). First,
  we show that~\ref{itm:MatCond_accr_conductivity} is satisfied by the
  following calculation,
  \begin{align}
    \re z \symbfcal M(z)&=\re {( z\symbfcal{M}_1(z)+ \symbfcal{M}_0(z))}\notag\\ &=\varepsilon_{\omega}\nu+\paran{\frac{\nu(\nu+b)a}{t^2+(\nu+b)^2}+\frac{a t^2}{t^2+(\nu+b)^2}}\notag\\
                        &=\underbrace{\frac{\nu^2a+a t^2}{t^2+(\nu+b)^2}}_{>0}+\underbrace{\frac{b\nu a}{t^2+(\nu+b)^2}+\varepsilon_{\omega}\nu}_{\Big\{\substack{
\,>0\text{ if } \nu >0\\\,<0\text{ if } \nu <0}}\,.\label{eq:deb-rezmz}
  \end{align}
  We see that $\re z \symbfcal M(z)$ is positive on $\R_{\geq 0}\times \R$.
  Next, we show that~\ref{itm:MatCond_conductivity} is fulfilled. We observe that
\[
\lim_{z \to 0} z \chi(z) = \lim_{z \to 0} \frac{z a}{z + b} = 0.
\]
Note that $\dom(\symbfcal M)=\C\setminus\{-b +0\iu\}$ and therefore
$\C_{\re>-\nu_{0}}\setminus \dom(\symbfcal{M}_0)=\{-b +0\iu\}$ for all
$\nu_{0}>0$ and
\[\dom(\symbfcal M)\cap\C_{\re >-\nu_0}=\C_{\re > -\nu_0}\setminus\{-b +0\iu\}\,.\]
Moreover, since the denominator $z + b$ does not vanish in $\C_{\re > -b}$, the
function $\chi(z)$ is analytic and bounded on $\C_{\re > -b}$.
Finally, for all $z \in \C_{\re > 0}$, we have
\[
\varepsilon_{\omega} + \re \chi(z) = \varepsilon_{\omega} + \re\left( \frac{a}{z + b} \right) \geq \varepsilon_{\omega}.
\]
This inequality holds because $\re(z + b) > b > 0$, ensuring that $\re\left( \frac{a}{z + b} \right) \geq 0$.

Therefore, all the conditions specified in~Remark~\ref{rem:maccret2} are satisfied by the Debye material law. By Theorems 3.10 and 3.15 in~\cite{dohnalWellPosednessExponentialStability2023}, we conclude that the Debye model is exponentially stable.
\end{proofp}

\begin{corollary}{Exponential stability of the Lorentz model}{exp-stab-lorentz}
Consider Problem~\ref{problem:evolution-disp} and suppose $\bar{\sigma}$,
$L_{1}=\emptyset$, $L_{2}=\brc{0}$ and $d_{0}=0$. Equation~\eqref{eq:disp-mod}
is equivalent to a Lorentz model, i.\,e.
\begin{equation}
  \label{eq:disp-mod-lor}
  \varepsilon(z) =  \varepsilon_{\omega}
  + \frac{c}{e+z f+z^{2}}\,.
\end{equation}
The associated material law is not exponentially stable according to Definition~\ref{def:exp-stab}.
\end{corollary}
The proof can be found in~\autocite[Appendix A.1]{dohnalWellPosednessExponentialStability2023}.

\paragraph*{Analytic correction to the Lorentz dispersion model}
\label{sec:orge475eec}
The Maxwell system combined with the Lorentz susceptibility fails to meet the
criteria for exponential stability. However, these criteria are broad and ensure
exponential decay of the solution for right-hand sides \(\symbf{\Phi}\) and
\(\symbf{\Psi}\) in \(L^2_{-\nu}\). When the Fourier–Laplace transform of the
right-hand side is centered around a specific \emph{frequency} \(z = z_0 \in \C\),
the exact configuration of the solution operator becomes inconsequential as
\(\abs{z}\) increases. In~\autocite[Appendix
A.2]{dohnalWellPosednessExponentialStability2023} this is addressed by
setting \(r \gg 1\) and considering the following modified material law for
\(z_0 = 0\),
\begin{equation}
  \symbfcal{M}_r(z) \coloneq \frac{c}{e + z^2 - f z} \left( 1 + \frac{z}{r} \right),
  \qquad \text{with} \quad e,\, f > 0.
  \label{eq:DL_modified}
\end{equation}
This modified material law \(\symbfcal{M}_r\) is bounded within a half-plane
$\C_{\re > -\nu_0}$ for some $\nu_0>0$. As shown in~\autocite[Lemma
A.1]{dohnalWellPosednessExponentialStability2023}, \(\symbfcal{M}_r\) satisfies
the condition~\ref{itm:MatCond_accr1} when
\(r > \frac{\omega_0^2}{2 \gamma}\) is sufficiently large. The second
condition~\ref{itm:MatCond_accr2} is also met, ensuring the exponential
stability of \(\symbfcal{M}_r\).
Additionally, Ionescu-Tira and Dohnal propose an analytic correction localized
around \(z_0\),
\begin{equation*}
  \symbfcal{M}_r(z) = \varepsilon_0 + \sum_{j=1}^n \frac{c_j}{e_j + z^2 - f_j z} \left( 1 + \frac{z - z_0}{r} \right).
\end{equation*}

Although the interaction between the parameters of the material law and \(z_0\)
is more complex, the conditions~\ref{itm:MatCond_accr1}
and~\ref{itm:MatCond_accr2} hold for small
\(c\)~\autocite[Appendix A.2]{dohnalWellPosednessExponentialStability2023}.

\paragraph*{Exponential stability of the generalized dispersion model}
The generalized dispersion model~\eqref{eq:disp-mod} is constructed as a linear
combination of a conductive term, Lorentz models, and Debye models. When the
modified Lorentz model~\eqref{eq:DL_modified} is employed, each of these
components is proven to be exponentially stable. Due to the linearity of the
overall system, this stability property extends naturally to the entire
generalized model. Consequently, we conclude that the generalized dispersion
model is exponentially stable. Note that, although
Problem~\ref{problem:evolution-disp} is restricted to isotropic material, an
extension to anisotropic material is straightforward. This extension lies
outside the scope of this work. Instead, we proceed by applying
Theorem~\ref{thm:picards} to PML and discuss the coupling of the physical and
the artificial PML problem.

\section{Perfectly matched layers \label{sec:pml-picard}}
\label{sec:orgd1a064f}
In numerical simulations, wave propagation must be truncated to bounded regions
to avoid unphysical reflections caused by boundary conditions. The Perfectly
Matched Layer (PML) is an artificial absorbing layer used to truncate unbounded
domains to bounded domains in wave propagation problems. By using an analytic
continuation of the solution to the complex plane, PMLs ensure that no
reflections occur and achieve exponential damping through an appropriate
coordinate transformation. We introduce a complex coordinate stretching
transformation that attenuates outgoing waves:
\[
\tilde{x} = \int_{0}^{x} s_x(\xi) \, \drv\xi, \quad \tilde{y} = \int_{0}^{y}
s_y(\eta) \, \drv\eta, \quad \tilde{z} = \int_{0}^{z} s_z(\zeta) \, \drv\zeta,
\]
where \(s_x(\xi)\), \(s_y(\eta)\), and \(s_z(\zeta)\) are complex-valued
stretching functions.

\paragraph*{Physical and PML domain} In the following, we decompose the domain
\(\Omega\) into two regions: the physical domain \(\Omega_{\text{phys}}\) and
the PML domain \(\Omega_{\text{PML}}\), i.e.,
\(\Omega = \Omega_{\text{phys}} \cup \partial\Omega_{\text{phys}} \cup
\Omega_{\text{PML}}\) with
\(\Omega_{\text{phys}} \cap \Omega_{\text{PML}} = \emptyset\). To ensure the
necessary regularity, we assume that \(\Omega_{\text{phys}}\) and
\(\Omega_{\text{PML}}\) are open, disjoint, and their union is the full domain
\(\Omega\). Moreover, the interface $\partial\Omega_{\text{phys}}$ between
\(\Omega_{\text{phys}}\) and \(\Omega_{\text{PML}}\) is sufficiently
Lipschitz-continuous.

\paragraph*{Complex frequency shifted perfectly matched layer (CFS-PML)}
\label{sec:orga8feaa5}
CFS-PMLs employ a frequency-shifted complex coordinate stretching function:
\begin{equation}\label{eq:cfspml1-a}
s_{k}(\nu) = 1 + \frac{\sigma_{k}}{\alpha_{k} + \iu\nu}, \quad k \in
\{x,\,y,\,z\}, \quad \vct s(\nu) =(s_{x}(\nu),\,s_{y}(\nu),\,s_{z}(\nu))
\end{equation}
The loss rate \(\sigma_k\) is chosen such that \(\sigma_{k}=0\) at the interface
to the physical domain, to ensure continuity of the physical parameters and
solutions. The frequency shift \(\alpha_k\) typically decays towards the PML
boundary. Again, the stretching function $s_{k}$ is extended into the physical
domain as the identity function. For brevity, the spatial dependence of \(s_k\),
\(\sigma_k\), and \(\alpha_k\) is often omitted.

In Definition~\ref{def:mw-operator}, we introduced a differential operator for
the Maxwell equation under coordinate transformations, which we now consider
under coordinate transformations. The physical coordinates and the complex
stretched coordinates are related by
\begin{equation}
  \label{eq:complex-coord-upml}
  \frac{\partial \tilde{x_k}}{\partial x_k} = s_k, \quad \frac{\partial \bullet}{\partial \tilde{x_k}} = \frac{1}{s_k} \frac{\partial \bullet}{\partial x_k}.
\end{equation}
\begin{definition}{Differential operators under coordinate transformations}{differential-operators-coord-curl}
  Let \(\Omega \subseteq \R^d\) be an open set. Under the coordinate transformations~\eqref{eq:complex-coord-upml}, we define:
  \begin{equation}
    \label{eq:curl-s}
    \curl_{\vct s}\colon \mat H(\curl,\,\Omega)^3 \subset L^2(\Omega)^3 \to L^2(\Omega)^3, \quad
    \symbf\varphi \mapsto \nabla \times \symbf\varphi = \begin{pmatrix}
      s_y^{-1} \partial_y \varphi_{3} - s_z^{-1} \partial_z \varphi_y \\
      s_z^{-1} \partial_z \varphi_{1} - s_x^{-1} \partial_x \varphi_z \\
      s_x^{-1} \partial_x \varphi_{2} - s_y^{-1} \partial_y \varphi_x
    \end{pmatrix}.
  \end{equation}
  The Jacobian matrix of the coordinate transformation is:
  \begin{equation}
  \label{eq:coord-jac}
  \mat{J} \coloneq s_{k} \vct e_k \otimes \vct e_k = \begin{pmatrix}
    s_{x} & 0 & 0 \\
    0 & s_{y} & 0 \\
    0 & 0 & s_{z}
  \end{pmatrix}.
  \end{equation}
  The coordinate transformations are expressed as:
  \begin{equation}
    \label{eq:jac-trafo-curl}
    \curl_{0,\,\vct s} \vct{E} = \operatorname{det} (\mat{J})^{-1} \mat{J} \curl \vct{E}.
  \end{equation}
\end{definition}

Considering the general Cauchy problem~\eqref{eq:Maxwell_pos_init} with a PML and using~\eqref{eq:jac-trafo-curl},
\begin{equation}
  \left( \partial_t \begin{pmatrix} \operatorname{det} (\mat{J}) \mat{J}^{-1}
                      \mat{J}^{-\top} \varepsilon(\partial_t) & 0 \\ 0 &
                                                                         \operatorname{det} (\mat{J}) \mat{J}^{-1} \mat{J}^{-\top} \mu \end{pmatrix} + \begin{pmatrix} 0 & \curl \\ \curl_0 & 0 \end{pmatrix} \right) \begin{pmatrix} \vct{E} \\ \vct{H} \end{pmatrix} = \begin{pmatrix} \vct f \\ 0 \end{pmatrix}.
  \label{eq:Mw1ord_2-pml1}
\end{equation}

For an isotropic PML in Cartesian coordinates, the material law \(\symbfcal{M}\)
is multiplied by $ \vct s$ component-wise, as \(\operatorname{det} (\mat{J})
\mat{J}^{-1} \mat{J}^{-\top} = \vct s\), where \(\vct s\) is defined as
in~\eqref{eq:cfspml1-a} for the CFS-PML, respectively,
\begin{equation}
  \left( \partial_t \begin{pmatrix} \vct s \varepsilon(\partial_t) & 0 \\ 0 & \vct s
                                                                         \mu \end{pmatrix}
                                                                       + \begin{pmatrix}
                                                                           0 &
                                                                               \curl \\ \curl_0 & 0 \end{pmatrix} \right) \begin{pmatrix} \vct{E} \\ \vct{H} \end{pmatrix} = \begin{pmatrix} \vct f \\ 0 \end{pmatrix}.
  \label{eq:Mw1ord_2-pml2}
\end{equation}

\section{A dispersive Maxwell equation with PML\label{sec:dispersion-pml-picard}}
\label{sec:orgded6bbd}
We formulate the systems for the CFS-PML, prove the well-posedness of
the arising problem and investigate its exponential stability. We then integrate the
dispersion models we investigated in Section~\ref{sec:dispersion-picard} into
the PML formulations.
\subsection{The CFS-PML for a non-dispersive Maxwell equation}
\label{sec:orgca87162}
To provide a foundational context and for purposes of completeness, we consider
the CFS-PML applied to the non-dispersive Maxwell equations. This serves as an
introductory example before addressing more complex models. In the CFS-PML
region we apply the coordinate transformation~\eqref{eq:cfspml1-a} to
Problem~\ref{problem:evolution-disp}, assuming $\bar{\sigma}=0$,
$L_{1}=\emptyset$, $L_{2}=\emptyset$. To solve the PDE numerically later on and
avoid evaluations of convolution terms, we need to introduce auxiliary
variables, due to the more involved stretching function.
\begin{problem}{Maxwell equation with CFS-PML}{pml-waves-1d}
  Let \(\symbfcal{H}\) denote the Hilbert space
  \(\symbfcal{H}\coloneq L(\Omega)^{4}\).
  Let $\symbfcal M_0,\: \symbfcal M_1\colon \symbfcal{H}\to \symbfcal{H}$ and
  $\symbfcal A\colon \dom(\symbfcal A)\subset \symbfcal{H} \to \symbfcal{H}$
  defined by
  \begin{equation}
    \label{eq:pml-evo}
    \symbfcal{M}_0=\begin{pmatrix}
            1 & {\color{gray} \symbf 0}      &{\color{gray} \symbf 0} &{\color{gray} \symbf 0} \\
            {\color{gray} \symbf 0}       &1 &{\color{gray} \symbf 0} &{\color{gray} \symbf 0} \\
            {\color{gray} \symbf 0}       &{\color{gray} \symbf 0}       &1&{\color{gray} \symbf 0} \\
            {\color{gray} \symbf 0}       &{\color{gray} \symbf 0}       &{\color{gray} \symbf 0} &1\\
          \end{pmatrix}\;,\;
          \symbfcal{M}_1=\begin{pmatrix}
                  \sigma  &{\color{gray} \symbf 0}        &-\alpha &{\color{gray} \symbf 0} \\
                  {\color{gray} \symbf 0}        &\sigma  &{\color{gray} \symbf 0}        &-\alpha \\
                  -\sigma &{\color{gray} \symbf 0}        &\alpha  &{\color{gray} \symbf 0} \\
                  {\color{gray} \symbf 0}        &-\sigma &{\color{gray} \symbf 0}        &\alpha \\
                \end{pmatrix}\;,\;
                \symbfcal A=\begin{pmatrix}
                    {\color{gray} \symbf 0}& -\curl &{\color{gray} \symbf 0} &{\color{gray} \symbf 0} \\
                    \curl_0 &{\color{gray} \symbf 0}      &{\color{gray} \symbf 0} &{\color{gray} \symbf 0} \\
                    {\color{gray} \symbf 0}&{\color{gray} \symbf 0}      &{\color{gray} \symbf 0} &{\color{gray} \symbf 0} \\
                    {\color{gray} \symbf 0}&{\color{gray} \symbf 0}      &{\color{gray} \symbf 0} &{\color{gray} \symbf 0} \\
    \end{pmatrix}\,,
  \end{equation}
  with the unknowns
  \begin{equation}
    \label{eq:pml-var}
    \vct U={(\vct E,\:\vct H,\:\vct R,\:\vct Q)}^{\top},\; \vct f={(0,\:0,\:0,\:0)}^{\top}\,.
  \end{equation}
  \emph{For a given right-hand side $\vct f\in L^2_{\nu}(\R;\:\symbfcal{H})$ find
    $\vct U\in L^2_{\nu}(\R;\:\symbfcal{H})$ such that}
  \begin{equation}
    \label{eq:evo-problem-cfs1d}
    (\partial_t \symbfcal M_0 + \symbfcal M_1 + \symbfcal A)\vct U=\vct f.
  \end{equation}
\end{problem}
\begin{corollary}{Well-posedness of Problem~\ref{problem:pml-waves-1d}}{well-posed-evo-c}
Problem~\ref{problem:pml-waves-1d} is well-posed. That is for each $\vct f\in L^2_{\nu}(\R;\:\symbfcal{H})$ there exists a unique solution $\vct U\in L^2_{\nu}(\R;\:\symbfcal{H})$.
\end{corollary}
\begin{proofp}
  The operator $\symbfcal A$ is skew self-adjoint, and $\symbfcal M_0$,
  $\symbfcal M_1$ are bounded linear operators with $\symbfcal M_0$ self-adjoint
  and strictly positive definite. Then it holds
  \begin{align*}
      \re \dd{z \symbfcal M(z) \vct u}{\vct u} &= \re \dd{z \symbfcal M_0 \vct u}{\vct u} + \dd{\re \symbfcal M_1 \vct u}{\vct u} \geq \nu \dd{\vct u}{\vct u} + \dd{\symbfcal M_1 \vct u}{\vct u} \\
      &\geq (\sigma - \frac{1}{2}\alpha - \frac{1}{2}\sigma + \nu) \dd{\vct E}{\vct E} + (\alpha - \frac{1}{2}\alpha - \frac{1}{2}\sigma + \nu) \dd{\vct R}{\vct R} \\
      &\quad + (\sigma - \frac{1}{2}\alpha - \frac{1}{2}\sigma + \nu) \dd{\vct H}{\vct H} + (\alpha - \frac{1}{2}\alpha - \frac{1}{2}\sigma + \nu) \dd{\vct Q}{\vct Q}.
  \end{align*}
  Clearly, there exists a $\nu_0 > 0$ such that~\eqref{eq:evo-cond} is fulfilled. Thus, Theorem~\ref{thm:picards} concludes the proof.
\end{proofp}
\paragraph*{Exponential stability of Maxwell equation with CFS-PML}
Using~\eqref{eq:no-accretivity}, we immediately see that exponential stability for the CFS-PML in the ADE formulation is not
possible. Instead of the formulation with ADEs in Problem~\ref{problem:pml-waves-1d}, we
consider the abstract material law in the Fourier-Laplace domain.
\begin{corollary}{Exponential stability of
    Problem~\ref{problem:pml-waves-1d}}{exp-stab-pml-prob}
  The CFS-PML, as stated in Problem~\ref{problem:pml-waves-1d} is exponentially
  stable.
\end{corollary}
The proof follows analogously to Corollary~\ref{cor:exp-stab-deb}. Note that
for the CFS-PML, \ref{itm:MatCond_accr_conductivity} and
\ref{itm:MatCond_conductivity} hold for both components of the material law
\(\vct e\) and \(\vct h\).

\subsection{The CFS-PML for a dispersive Maxwell equation}
\label{sec:org2614122}
In this section, we apply the coordinate
transformation~\eqref{eq:cfspml1-a} in the CFS-PML region to Problem~\ref{problem:evolution-disp}. According to
Definition~\ref{def:differential-operators-coord-curl} this yields the system of equations
\begin{equation}
  \label{eq:disp-pde-pml}
  \begin{aligned}
     \paran{\varepsilon_{\omega} \partial_{t}+\sigma \varepsilon_{\omega} +\bar{\sigma} + \sum_{l\in L_{2}}d_{l}}\vct E +
     \sum_{l\in L_{1}} (a_l\vct E-(b_l-\sigma)\vct p_l) +
     \sum_{l\in L_{2}} (\vct j_{l} +\sigma \vct p_l)+
    \vct S_{1} &-\alpha(\vct S_2+\vct S_3) -\curl \vct H=\vct f\,,\\
    \partial_{t}\vct H - \alpha \vct R + \sigma \vct H + \curl_0 \vct E &= 0\,,\\
    \partial_t \vct p_{l} - a_{l} \vct E + b_{l} \vct p_{l} &= 0\,,\; l\in L_{1}\,,\\
    \partial_{t} \vct j_m + (d_mf_m - c_m)\vct E + f_m\vct j_m + e_m\vct p_m=0\,,\quad\partial_{t} \vct p_m - d_m \vct E - \vct j_m &= 0\,,\; m\in L_{2}\,,\\
    \partial_t \vct S_{1} + \alpha \vct S_{1} -\sigma\bar{\sigma}\vct E&=0\,,\\
    \partial_t \vct S_{2} + \alpha \vct S_{2} -\varepsilon_{\omega}\sigma \vct E&=0\,,\\
    \partial_t \vct S_{3} + \alpha \vct S_{3} -\sigma \sum_{l\in L_1\cap L_2}
    \vct p_l&=0\,,\\
    \partial_{t}\vct R+\alpha \vct R - \sigma \vct H &= 0\,.
  \end{aligned}
\end{equation}
We recast~\eqref{eq:disp-pde-pml} as an evolutionary problem by putting
{
\begin{subequations}\label{eq:disp-ops-pml}\allowdisplaybreaks
  \begin{equation}
    \symbfcal A=\begin{tikzpicture}[baseline=(current bounding box.center)]%
      \matrix[sparsemat] {
         &-\curl &\cdots & \\
        \curl_0& &\cdots  &  \\
        \vdots&\vdots   & \ddots & \vdots \\
        &   &\cdots &   \\
      };
    \end{tikzpicture}\;,\quad
    \symbfcal{M}_0=\begin{tikzpicture}[baseline=(current bounding box.center)]%
      \matrix[sparsemat] {
        \varepsilon_{\omega} & &\cdots & \\
        & 1 &\cdots  &  \\
        \vdots&\vdots   & \ddots & \vdots \\
        &   &\cdots &   1\\
      };
    \end{tikzpicture}
    \normalsize
  \end{equation}
  \begin{equation}
    \symbfcal{M}_1=\begin{tikzpicture}[baseline=(current bounding
      box.center),font=\footnotesize]%
      \matrix[sparsemat] {
        (\bar{\sigma} + \rho)+\varepsilon_{\omega}\sigma &         & -\vec(b_{L_1}-\sigma)^{\top}  &               1&  \vec(\sigma)^{\top} & 1      & -\alpha & -\alpha& \\
                                                        & \sigma  &                               &                &                      &        &         &        &-\alpha \\
                                       -\vec(a_{L_{1}})  &         & \diag(b_{L_1})                &                &                      &        &         &        & \\
                      \vec(d_{L_2}f_{L_2})-\vec(c_{L_2})  &         &                               & \diag(f_{L_2}) & \diag(e_{L_2})       &        &         &        & \\
                                         -\vec(d_{L_2})  &         &                               & -1             &                      &        &         &        & \\
                                    -\sigma\bar{\sigma}  &         &                               &                &                      & \alpha &         &        & \\
                             -\sigma\varepsilon_{\omega}  &         &                               &                &                      &        &  \alpha &        & \\
                                    -\sigma\bar{\sigma}  &         & -\vec(\sigma)^{\top}          &                & -\vec(\sigma)^{\top} &        &         & \alpha & \\
                                                         & -\sigma &                               &                &                      &        &         &        & \alpha \\
      };
    \end{tikzpicture}
    \normalsize
  \end{equation}
\end{subequations}
}
and
\begin{equation}
  \label{eq:disp-sols-pml}
  \vct U={(\vct E,\,\vct H,\,\vct j_{L_{1}},\,\,\vct j_{L_{2}},\,\vct
    p_{L_{2}},\,\vct S_1,\,\vct S_2\,\vct S_3,\,\vct R)}^{\top},\; \vct f={(0,\,0,\,0,\,0,\,0,\,0,\,0,\,0,\,0)}^{\top}.
\end{equation}
\begin{problem}{Dispersive Maxwell equation with PML as an evolutionary
    problem}{evolution-disp-pml}
  Let \(\symbfcal{H}\) denote the Hilbert space
  \(\symbfcal{H}\coloneq L^{9d}\),
  equipped with the inner product of $L^{9d}$. Let
  $\symbfcal M_0,\: \symbfcal M_1\colon \symbfcal{H}\to \symbfcal{H}$ and $\symbfcal A\colon \dom(\symbfcal A)\subset \symbfcal{H} \to \symbfcal{H}$ be defined by~\eqref{eq:disp-ops-pml}. The
  evolutionary problem with $\vct f$ and $\vct U$ defined according
  to~\eqref{eq:disp-sols-pml} reads:
  \emph{For a given right-hand side $\vct f\in L^2_{\nu}(\R;\:\symbfcal{H})$ find
    $\vct U\in L^2_{\nu}(\R;\:\symbfcal{H})$ such that}
  \begin{equation}
    \label{eq:evo-problem-disppml}
    (\partial_t \symbfcal M_0 + \symbfcal M_1+ \symbfcal A)\vct U=\vct f.
  \end{equation}
\end{problem}
\begin{corollary}{Well-posedness of
    Problem~\ref{problem:evolution-disp-pml}}{well-posed-evo-disp-pml}
  Problem~\ref{problem:evolution-disp-pml} is well-posed. That is for each
  $\vct f\in L^2_{\nu}(\R;\:\symbfcal{H})$ there exists a unique solution $\vct U\in L^2_{\nu}(\R;\:\symbfcal{H})$.
\end{corollary}
\begin{proofp}
  For Problem~\ref{problem:evolution-disp-pml} operator $\symbfcal A$ is skew
  self-adjoint. The operators $\symbfcal M_0$, $\symbfcal M_1$ are clearly bounded linear operators
  and $\symbfcal M_0$ is self-adjoint and positive definite. Applying the same
  technique as for the proof of Corollary~\ref{cor:well-posed-evo-c}, we are able to
  find a $\nu_0>0$ and $\gamma>0$ such that~\eqref{eq:evo-cond} holds. Therefore,
  Theorem~\ref{thm:picards} concludes the proof.
\end{proofp}
We establish exponential stability for the Lorentz and Debye models,
which are included in the formulation of Problem~\ref{problem:evolution-disp-pml}.
\begin{corollary}{Exponential stability of the Debye and Lorentz models with CFS-PML}{exp-stab-lor-deb-pml}
  Consider Problem~\ref{problem:evolution-disp-pml} and suppose either
\begin{equation}
  \bar{\sigma}=0,\quad L_{1}=\brc{0},\:L_{2}=\emptyset\text{ and }d_{0}=0\,,\label{eq:dm-lor}
\end{equation}
or
\begin{equation}
  \bar{\sigma}=0,\quad L_{1}=\brc{0},\:L_{2}=\emptyset\,.\label{eq:dm-deb}
\end{equation}
  Then, with~\eqref{eq:dm-lor}, Equation~\eqref{eq:disp-mod} is equivalent to a
  Lorentz model and the CFS-PML with Lorentz model is exponentially stable.
  Otherwise, with~\eqref{eq:dm-deb}, Equation~\eqref{eq:disp-mod} is equivalent to a
  Debye model and the CFS-PML with Debye model is exponentially stable.
\end{corollary}
\begin{proofp}
  We consider the abstract material law resulting from the multiplication of the
  PML function with a material law of the form $\symbfcal{M}(z)=\varepsilon_{\omega}+ \chi(z)$. Then multiplication by the PML function yields
  \begin{equation*}
  \symbfcal M_{\text{PML}}(z)=z(1+\frac{\sigma}{\alpha + z})(\varepsilon_{\omega}+ \chi(z))\quad z\in \C\,.
  \end{equation*}
  We show
  that~\ref{pro:exp-stab-pos-def} is satisfied.
  \begin{align}
    \re {(\nu + \iu t)\symbfcal M_{\text{PML}}(\nu + \iu
    t)}&=\underbrace{\paran{\frac{\sigma t^2+\sigma
         \nu^2}{t^2+(\nu+\alpha)^2}+\frac{\sigma\alpha\nu}{t^2+(\nu+\alpha)^2}+\nu}}_{\customlabel{eq:ub-a}{(\text{A})}}\underbrace{\paran{\varepsilon_{\omega}+\re{\chi(\nu+\iu
         t)}}}_{\customlabel{eq:ub-b}{(\text{B})}}\notag\\
&\quad \underbrace{- \paran{1+\frac{\alpha\sigma}{t^2+(\nu+\alpha)^2}}t\im
         {\chi(\nu+\iu
  t)}}_{\customlabel{eq:ub-c}{(\text{C})}}\label{eq:lor-cfspml-accr1}
  \end{align}
  The arguments to show the exponential stability of the material law
  $\symbfcal M_{\text{PML}}$ are similar to the ones we used in
  Corollary~\ref{cor:exp-stab-lorentz-upml}. We consider the three
  terms~\ref{eq:ub-a},~\ref{eq:ub-b} and~\ref{eq:ub-c} separately.
  \paragraph*{Term~\ref{eq:ub-a}} We establish that there exists a $\tilde\nu_0$ such that term~\ref{eq:ub-a} is
  strictly positive on $\C_{\re> -\tilde\nu_0}$. 
  \begin{equation}
    \label{eq:est-ub-a}
\underbrace{\frac{\sigma t^2+\sigma
         \nu^2}{t^2+(\nu+\alpha)^2}}_{>0}+\underbrace{\frac{\sigma\alpha\nu}{t^2+(\nu+\alpha)^2}+\nu}_{\Big\{\substack{
\,>0\text{ if } \nu >0\\\,<0\text{ if } \nu <0}}
  \end{equation}
  The first term is always positive and~\ref{eq:ub-a} is continuous on
  $\C_{\re> -a}$. Thus, we find that there exist a $\tilde \nu_0$, such that the
  term~\ref{eq:ub-a} is positive on $\C_{\re >-\nu_0}$.
  \paragraph*{Term~\ref{eq:ub-b}}
  We establish that there exists a $\bar\nu_0$ such that term~\ref{eq:ub-b} is
  strictly positive on $\C_{\re> -\bar\nu_0}$. We distinguish the cases when
  $\chi$ represents a Debye or a Lorentz model.
  \subparagraph*{Case 1: Lorentz model}
  \begin{equation}
    \label{eq:lorentz-ub-b}
    \varepsilon_{\omega} + \re
    \chi(\nu + \iu t)> 0\qquad\text{if } \re
    \chi(\nu + \iu t) >
    -\varepsilon_{\omega}\,.
  \end{equation}
  \subparagraph*{Case 2: Debye model}
  \begin{equation}
    \label{eq:deb-ub-b}
    \varepsilon_{\omega} + \re
    \chi(\nu + \iu t)> 0\qquad\text{if } \nu >
    -a
  \end{equation}
  Note that for both models $\chi(z)$ has poles, which lie in the half-plane
  $\C_{\re<0}$ and $\chi(z)$ is analytic and bounded on $\C_{\re \geq 0}$. Then,
  by continuity, $\varepsilon_{\omega}+\re \chi(\nu + \iu t)$ remains positive
  for $\nu < 0$ large enough.
  
  Since $\varepsilon_{\omega}>0$, there exists a $\bar \nu_0$ such that the
  product of~\ref{eq:ub-a} and~\ref{eq:ub-b} is positive and bounded by a
  constant, if $\chi$ represents a Lorentz or a Debye model.
  There, we used the positivity
  of $\re \chi(\nu + \iu t)$ on $\R_{\geq 0}\times \R$ and the continuity.

  \paragraph*{Term~\ref{eq:ub-c}}
  We establish, that there exists a $\hat\nu_0$ such that term~\ref{eq:ub-c} is
  strictly positive on $\C_{\re> -\hat\nu_0}$. For this, we separate the
  cases of Debye and Lorentz models:
  \subparagraph*{Case 1: Lorentz model}
  \begin{equation}
    \label{eq:lorentz-ub-c}
    -\Bigg(\underbrace{1+\frac{\alpha\sigma}{t^2+(\nu+\alpha)^2}}_{\geq
      1}\Bigg)t\im \chi(\nu + \iu t) \geq -t\im \chi(\nu + \iu
    t)\underset{\\\text{ if } \nu >
    -\frac{f}{2}}{>}
    0
  \end{equation}
  \subparagraph*{Case 2: Debye model}
  \begin{equation}
    \label{eq:deb-ub-c}
      -\paran{1+\frac{\alpha\sigma}{t^2+(\nu+\alpha)^2}}t\im \chi(\nu + \iu t)\geq -t\im \chi(\nu + \iu t)\geq 0
  \end{equation}
  By~\eqref{eq:lorentz-ub-c} and~\eqref{eq:deb-ub-c} we conclude that
  Proposition~\ref{pro:exp-stab-pos-def} is fulfilled; choosing a $\nu_0$,
  $0<\nu_0<\min\{\hat \nu_0,\,\tilde \nu_0,\,\bar \nu_0\}$, we find a $c_0$ such
  that
  \begin{equation}
    \label{}
  \re {(\nu + \iu t)\symbfcal M_{\text{PML}}(\nu + \iu t)}
                                 \ge c_0\text{ for all }\nu > -\nu_0\,.\label{eq:cfs-cfspml-exp-stab-lor}
  \end{equation}

  The second variable requires no treatment as the exponential stability can be
  shown analogously to~\eqref{eq:deb-rezmz}. Consequently, the Lorentz and Debye
  models combined with the CFS-PML are exponentially stable.
\end{proofp}

\section{Conclusion}
We presented a Hilbert space framework for analyzing dispersive Maxwell
equations, including the Debye and Lorentz models. We demonstrated the
well-posedness of these equations in an abstract setting and showed that the
resulting systems are both well-posed and exponentially stable, incorporating
perfectly matched layers (PMLs), specifically the complex frequency-shifted PML
(CFS-PML).

Our study addressed nonlinear perturbations and saturable nonlinearities in
electromagnetics. The results extend naturally to anisotropic or inhomogeneous
media. Careful treatment of material laws and detailed examination of
accretivity conditions were crucial in establishing exponential stability,
particularly when integrating PML into dispersive media. The results provide a
mathematical foundation for future research on numerical methods involving
dispersive Maxwell equations with perfectly matched layers. These methods have
the potential to significantly impact physics research in the design of novel
optical sources, optical waveguides, and studies in electromagnetic
compatibility.

\paragraph*{Competing interests}
The authors declare that there are no conflicts of interest related to this article.
\printbibliography
\appendix
\section{Uniaxial perfectly matched layer (UPML)}
\label{sec:org8137c3e}
Although CFS-PMLs are effective for handling anisotropic and dispersive media,
UPMLs offer a simpler formulation and are particularly attractive for their
computational efficiency and ease of implementation. These properties make UPMLs
a compelling alternative for some applications, especially when the objective is
to achieve a balance between accuracy and computational cost. UPMLs use a
specific form of the stretching function:
\begin{equation}\label{eq:upml1}
s_{k}(\nu) = 1 + \frac{\sigma_{k}}{\iu\nu}, \quad k \in \{x,\,y,\,z\},
\end{equation}
where the loss rate \(\sigma_k\) is a function of the position in the PML, where
\(\sigma_{k}=0\) at the interface
to the physical domain. The stretching function $s_{k}$ is extended into the physical domain as the identity function. UPMLs are widely used in computational electromagnetics due to their straightforward implementation and low memory requirements~\autocite{pledReviewRecentDevelopments2021}.

\subsection{The UPML for a non-dispersive Maxwell equation}
\label{sec:org0fed084}
In the UPML region, we apply the coordinate transformation~\eqref{eq:upml1} to
Problem~\ref{problem:evolution-disp}, with $\bar{\sigma}=0$, $L_{1}=\emptyset$,
and $L_{2}=\emptyset$. There is no source term in the PML region, i.e., \(\vct f=0\).
\begin{problem}{Maxwell equation with UPML}{upml-waves-1d}
  Let \(\symbfcal{H}\) denote the Hilbert space
  \(\symbfcal{H}\coloneq L(\Omega)^{6d}\).
  Let $\symbfcal M_0,\: \symbfcal M_1\colon \symbfcal{H}\to \symbfcal{H}$ be
  defined as
  \begin{equation}
    \label{eq:upml-evo}
    \symbfcal{M}_0=\begin{pmatrix}
            1 &{\color{gray} \symbf 0}       \\
                   {\color{gray} \symbf 0}&1 \\
          \end{pmatrix}\;,\quad
          \symbfcal{M}_1=\begin{pmatrix}
                  \sigma  &{\color{gray} \symbf 0}        \\
                         {\color{gray} \symbf 0} &\sigma  \\
                \end{pmatrix}\,,
  \end{equation}
  with the unknowns
  \begin{equation}
    \label{eq:upml-var-upml1d}
    \vct U=\begin{pmatrix}
    \vct E\\\vct H
    \end{pmatrix},\quad \vct f=\begin{pmatrix}
    \vct 0\\\vct 0
    \end{pmatrix}\,.
  \end{equation}
  \emph{For a given right-hand side $\vct f\in L^2_{\nu}(\R;\:\symbfcal{H})$ find
    $\vct U\in L^2_{\nu}(\R;\:\symbfcal{H})$ such that}
  \begin{equation}
    \label{eq:evo-problem-c}
    (\partial_t \symbfcal M_0 + \symbfcal M_1 + \symbfcal A)\vct U=\vct f.
  \end{equation}
\end{problem}
\begin{corollary}{Well-posedness of Problem~\ref{problem:upml-waves-1d}}{well-posed-evo-b}
Problem~\ref{problem:upml-waves-1d} is well-posed. That is for each $\vct f\in L^2_{\nu}(\R;\:\symbfcal{H})$ there exists a unique solution $\vct U\in L^2_{\nu}(\R;\:\symbfcal{H})$.
\end{corollary}
\begin{proofp}
  The operator $\symbfcal A$ is skew self-adjoint, and $\symbfcal M_0$, $\symbfcal M_1$ are bounded linear operators with $\symbfcal M_0$ self-adjoint and strictly positive definite, and $\symbfcal M_1$ positive definite. Since
  \begin{equation*}
      \dd{(\nu \symbfcal M_0 + \re \symbfcal M_1) \vct u}{\vct u}_{\symbfcal{H}} \geq \nu \dd{\vct u}{\vct u} + \dd{\re \symbfcal M_1 \vct u}{\vct u} \geq (\nu + \sigma) \dd{\vct u}{\vct u} \geq \nu_0 \dd{\vct u}{\vct u}\,,
  \end{equation*}
  condition~\eqref{eq:evo-cond} holds for $0 < \gamma < \nu_{0}$. Thus, Theorem~\ref{thm:picards} concludes the proof.
\end{proofp}
\begin{corollary}{Exponential stability of
    Problem~\ref{problem:upml-waves-1d}}{exp-stab-upml-prob}
  The Problem~\ref{problem:upml-waves-1d} is exponentially stable with decay rate
  $\sigma$ according to Definition~\ref{def:exp-stab}.
\end{corollary}
\begin{proofp}
  By Proposition~\ref{pro:exp-stab-pos-def}, $\dom(\symbfcal M) = \C$, and $\C_{\re > -\nu_{0}} \setminus \dom(\symbfcal{M}_0) = \{0\}$ for all $\nu_{0} > 0$. For $\nu > \nu_{0}$,
  \begin{align*}
    \label{stab-cond-ineq}
    \re \dd{z \symbfcal M(z) \vct u}{\vct u}_{\symbfcal{H}} &= \re \dd{z \symbfcal M_0 \vct u}{\vct u}_{\symbfcal{H}} + \dd{\symbfcal M_1 \vct u}{\vct u}_{\symbfcal{H}} \notag \\
    &> -\nu \dd{\vct u}{\vct u}_{\symbfcal{H}} + \sigma \dd{\vct u}{\vct u}_{\symbfcal{H}} = (\sigma - \nu) \dd{\vct u}{\vct u}_{\symbfcal{H}} \geq \left( \min\limits_{u \in \Omega} \sigma(\vct x) - \nu \right) \dd{\vct u}{\vct u}_{\symbfcal{H}}.
  \end{align*}
  Defining the decay rate as $\nu_0 = \min\limits_{\vct x \in \Omega}
  \sigma(\vct x) - c$ ensures that~\eqref{exp-stab-cond} holds.
\end{proofp}

\subsection{The UPML for a dispersive Maxwell equation}
\label{sec:org620a2d4}
We apply the coordinate
transformation~\eqref{eq:upml1} in the UPML region to the physical
problem~\ref{problem:evolution-disp} and which yields the equations
\begin{equation}
  \label{eq:disp-pde-upml}
  \begin{aligned}
    \varepsilon_{\omega} (\partial_{t} + \sigma ) \vct E +  \paran{\bar{\sigma} + \sum_{l\in L_{2}}d_{l}}\vct E +
    \sum_{l\in L_{1}} (a_l\vct E-b_l\vct p_l) +
     \sum_{l\in L_{2}} \vct j_{l}+
    \sigma \sum_{l\in L_1\cup L_{2}} \vct p_l&+
    \vct S -\curl \vct H=\vct f\,,\\
     \partial_{t}\vct H + \sigma \vct H + \Grad_0 \vct E &= 0\,,\\
    \partial_t \vct p_{l} - a_{l} \vct E + b_{l} \vct p_{l} &= 0\,,\quad l\in L_{1}\,,\\
    \partial_{t} \vct j_m + (d_mf_m - c_m)\vct E + f_m\vct j_m +
    e_m\vct p_m=0\,,\quad \partial_{t} \vct p_m - d_m \vct E - \vct j_m &= 0\,,\quad l\in L_{2}\,,\\
    \partial_t \vct S -\sigma\bar{\sigma}\vct E&=0\,.
  \end{aligned}
\end{equation}
We recast Equation~\eqref{eq:disp-pde-upml} as an evolutionary problem by
putting
{\footnotesize
\begin{subequations}\label{eq:disp-ops-upml}\allowdisplaybreaks
  \begin{equation}
    \symbfcal A=\begin{tikzpicture}[baseline=(current bounding box.center)]%
      \matrix[sparsemat] {
         &-\curl &\cdots & \\
        \curl_0& &\cdots  &  \\
        \vdots&\vdots   & \ddots & \vdots \\
        &   &\cdots &   \\
      };
    \end{tikzpicture}\;,\quad
    \symbfcal{M}_0=\begin{tikzpicture}[baseline=(current bounding box.center)]%
      \matrix[sparsemat] {
        \varepsilon_{\omega} & &\cdots & \\
        & 1 &\cdots  &  \\
        \vdots&\vdots   & \ddots & \vdots \\
        &   &\cdots &   1\\
      };
    \end{tikzpicture}
    \normalsize
  \end{equation}
  \begin{equation}
    \symbfcal{M}_1=\begin{tikzpicture}[baseline=(current bounding
      box.center)]%
      \matrix[sparsemat] {
         (\bar{\sigma} + \rho)+\varepsilon_{\omega}\sigma &         & - \vec(b_{L_1})^{\top}+\sigma  & 1& \sigma  & 1 \\
        & \sigma &         &              & &  \\
        -\vec(a_{L_{1}})  & & \diag(b_{L_1})       &  &   &  \\
        \vec(d_{L_2}f_{L_2})-\vec(c_{L_2}) & &         & \diag(f_{L_2}) & \diag(e_{L_2}) &\\
        -\vec(d_{L_2})  & & & -1& & \\
        -\sigma\bar{\sigma}  & &  & & & \\
      };
    \end{tikzpicture}
    \normalsize
  \end{equation}
\end{subequations}
}
and
\begin{equation}
  \label{eq:disp-sols-upml}
  \vct U={(\vct E,\,\vct H,\,\vct j_{L_{1}},\,\,\vct j_{L_{2}},\,\vct p_{L_{2}},\,\vct S)}^{\top},\; \vct f={(0,\,0,\,0,\,0,\,0,\,0)}^{\top}.
\end{equation}
\begin{problem}{Dispersive Maxwell equation with PML as an evolutionary
    problem}{evolution-disp-upml}
  Let \(\symbfcal{H}\) denote the Hilbert space
  \(\symbfcal{H}\coloneq L^{6d}\)
  equipped with the inner product of $L^{6d}$. Let
  $\symbfcal M_0,\: \symbfcal M_1\colon \symbfcal{H}\to \symbfcal{H}$ and $\symbfcal A\colon \dom(\symbfcal A)\subset \symbfcal{H} \to \symbfcal{H}$ be defined by~\eqref{eq:disp-ops-upml}. The
  evolutionary problem with $\vct f$ and $\vct U$ defined according
  to~\eqref{eq:disp-sols-upml} reads:
  \emph{For a given right-hand side $\vct f\in L^2_{\nu}(\R;\:\symbfcal{H})$ find
    $\vct U\in L^2_{\nu}(\R;\:\symbfcal{H})$ such that}
  \begin{equation}
    \label{eq:evo-problem-dispupml}
    (\partial_t \symbfcal M_0 + \symbfcal M_1+ \symbfcal A)\vct U=\vct f.
  \end{equation}
\end{problem}
\begin{corollary}{Well-posedness of
    Problem~\ref{problem:evolution-disp-upml}}{well-posed-evo-disp-upml}
  Problem~\ref{problem:evolution-disp-upml} is well-posed. That is for each
  $\vct f\in L^2_{\nu}(\R;\:\symbfcal{H})$ there exists a unique solution $\vct U\in L^2_{\nu}(\R;\:\symbfcal{H})$.
\end{corollary}
\begin{proofp}
  For Problem~\ref{problem:evolution-disp-upml} operator $\symbfcal A$ is skew
  self-adjoint. The operators $\symbfcal M_0$, $\symbfcal M_1$ are clearly bounded linear operators
  and $\symbfcal M_0$ is self-adjoint and positive definite. Applying the same
  technique as for the proof of Corollary~\ref{cor:well-posed-evo-c}, we are able to
  find a $\nu_0>0$ and $\gamma>0$ such that~\eqref{eq:evo-cond} holds. Therefore,
  Theorem~\ref{thm:picards} concludes the proof.
\end{proofp}
Similar to Corollary~\ref{cor:exp-stab-deb} and~\ref{cor:exp-stab-lorentz}, we
establish exponential stability for parameters corresponding to a Debye and a
Lorentz model in Problem~\ref{problem:evolution-disp-upml}.
\begin{corollary}{Exponential stability of the Lorentz and the Debye model with UPML}{exp-stab-lorentz-upml}
  Consider Problem~\ref{problem:evolution-disp-upml} and suppose $\bar{\sigma}=0$,
  $L_{1}=\emptyset$, $L_{2}=\brc{0}$ and $d_{0}=0$. Then,~\eqref{eq:disp-mod}
  is equivalent to a Lorentz model. Then Problem~\ref{problem:evolution-disp-upml}  is exponentially stable.
\end{corollary}
\begin{proofp}
  We consider the abstract material law resulting from the multiplication of the
  PML function with the Lorentz material law given by
  \begin{equation*}
  \symbfcal M_{\text{PML}}(z)=z(1+\frac{\sigma}{z})\symbfcal M(z)\quad z\in \C\,.
  \end{equation*}
  We show that Proposition~\ref{pro:exp-stab-pos-def} is satisfied by the following calculation,
  \begin{align}
    \re {(\nu + \iu t)\symbfcal M(\nu + \iu t)}&=\nu\re {\symbfcal M(\nu + \iu
                                                 t)-t\im \symbfcal M(\nu + \iu
                                                 t)+\sigma \re \symbfcal M(\nu +
                                                 \iu t)}\,.\\
    \intertext{Now we need to distinguish the cases if $\chi$ corresponds to a
    Lorentz or a Debye model.
    \paragraph*{Case 1: Lorentz Model} If $\chi$ corresponds to a Lorentz model, we have}
    \re {(\nu + \iu t)\symbfcal M(\nu + \iu t)}&=\underbrace{(\nu + \sigma)}_{\substack{> 0\\ \text{
    if } \nu > -\sigma}}%
    \underbrace{(\varepsilon_{\omega} + \re
    \chi(\nu + \iu t))}_{\substack{> 0\\ \text{ if } \re
    \chi(\nu + \iu t) >
    -\varepsilon_{\omega}}}%
    \underbrace{-t\im \chi(\nu + \iu t)}_{\substack{>
    0\\ \text{ if } \nu >
    -\frac{f}{2}}}\label{eq:lor-pml-accr1}\,.\\
    \intertext{Proposition~\ref{pro:exp-stab-pos-def} is fulfilled as the first
    and third terms are positive for some \(\nu\), \(0 < \nu < \min(\sigma,\frac{f}{2})\).
    For the middle term, note that $\re \chi(\nu + \iu t)$ is positive and bounded on $\R_{\geq 0}\times \R$.
    By continuity, $\varepsilon_{\omega}+\re \chi(\nu + \iu t)$ remains positive for $\nu < 0$ small enough.
    Thus, there exists \(\nu_0 \in (0, \min(\sigma, \frac{f}{2}))\) and \(c_0 > 0\) such
    that}
    \re {(\nu + \iu t)\symbfcal M_{\text{PML}}(\nu + \iu t)}
                                               &\ge c_0\text{ for all }\nu >
                                                 -\nu_0\,.\label{eq:cfs-pml-exp-stab-lor}
    \intertext{\paragraph{Case 2: Debye Model} If $\chi$ corresponds to a Debye model, we have}
    \re {(\nu + \iu t)\symbfcal M(\nu + \iu t)}&=\underbrace{(\nu + \sigma)}_{\substack{> 0\\ \text{
    if } \nu > -\sigma}}%
    (\varepsilon_{\omega} + \underbrace{\re
    \chi(\nu + \iu t)}_{\substack{> 0\\ \text{ if } \nu >
    -a}})%
    \underbrace{-t\im \chi(\nu + \iu t)}_{\geq
    0}\,.\label{eq:deb-pml-accr1}\\
    \intertext{Proposition~\ref{pro:exp-stab-pos-def} is fulfilled as all terms are positive for \(\nu > -\min(\sigma, a)\). The zeros of the denominators of \(\im \symbfcal M(\nu + \iu t)\) and \(\re \symbfcal M(\nu + \iu t)\) are at \(-a\). Thus, there exists \(\nu_0 \in (0, \min(\sigma, a))\) and \(c_0 > 0\) such that}
    \re {(\nu + \iu t)\symbfcal M_{\text{PML}}(\nu + \iu t)}
                                               &\ge c_0\text{ for all }\nu > -\nu_0\,.\label{eq:cfs-pml-exp-stab-deb}
  \end{align}
  The second variable requires no further treatment as the material law in this
  component is unaffected and analogous to~\ref{problem:upml-waves-1d}.
  Therefore, both Lorentz and Debye models combined with UPML are exponentially
  stable.
\end{proofp}
\end{document}